\documentclass[preprint,10pt]{elsarticle}

\usepackage{natbib}        % required for bibliography

\usepackage{amsfonts,amssymb,amstext,verbatim,amsmath,graphicx,color}
\usepackage{url}

\usepackage{amssymb}
\usepackage{latexsym}

\usepackage{graphics} %
\usepackage{epsfig} %
\usepackage{epstopdf} %

\newcommand{\ba}{\begin{array}}
\newcommand{\ea}{\end{array}}
\newcommand{\be}{\begin{equation}}
\newcommand{\ee}{\end{equation}}

\newtheorem{theorem}{\textbf{Theorem}}[section]
\newtheorem{lemma}[theorem]{\textbf{Lemma}}
\newtheorem{proposition}[theorem]{\textbf{Proposition}}

\newtheorem{assumption}[theorem]{\textbf{Assumption}}
\newtheorem{algo}[theorem]{\textbf{Algorithm}}

\begin{document}
\begin{frontmatter}
%\title{Using the Style File IEEEtran.sty}
\title{Hamiltonian-Based Algorithm for Optimal Control}


%\thanks[footnoteinfo]{Research of the first author was supported in
%part by  the NSF under Grant CNS-1239225.}

\tnotetext[]{{Research  supported in
		part by  the NSF under Grant CNS-1239225.}  }

\author[First]{M.T. Hale}
\author[First]{Y. Wardi}
\author[Second]{H. Jaleel}
\author[First]{M. Egerstedt}

\address[First]{School of Electrical and Computer Engineering,
Georgia Institute of Technology, Atlanta, GA 30332,
USA \\ Email:
\{mhale30,ywardi,magnus\}@ ece.gatech.edu.}
\address[Second]{Department of Electrical Engineering, University of Engineering and
Technology, Lahore, Pakistan\\
 Email: hassanjaleel@uet.edu.pk.}

%\maketitle
%\thispagestyle{plain}\pagestyle{plain}

\begin{abstract}
This paper proposes an algorithmic technique for a class of optimal control problems where it is easy to compute
a pointwise  minimizer of the Hamiltonian associated with every applied
control.
The  algorithm operates in the space of relaxed
controls and projects the final result into the space of ordinary controls. It is based on the descent direction from a given
relaxed control towards a pointwise minimizer of the Hamiltonian. This direction comprises a form of gradient projection and
for some systems,  is argued to have computational advantages over direct gradient directions. The algorithm is shown to be applicable to a class of hybrid optimal control problems.  The  theoretical results, concerning convergence of the algorithm,
are corroborated by simulation examples on switched-mode hybrid systems as well as on a problem of balancing transmission- and motion
energy in a mobile robotic system.

 \end{abstract}

\end{frontmatter}

\section{Introduction}

Consider  dynamical systems described by the differential equation
\begin{equation}
\dot{x}\ =\ f(x,u),
\end{equation}
where $x(t)\in R^n$ is its state variable and $u(t)\in U\subset R^k$ is the  input control.
The set $U\subset R^k$ is assumed to be compact.
Suppose that the initial time is $t_{0}=0$,  and the initial state $x_{0}:=x(0)\in R^n$ and
the final time $t_{f}>0$ are given and fixed. A control
$\{u(t),\  t\in[0,t_{f}]\}$ is said to be
{\it ordinary} if the function
$u:[0,t_{f}]\rightarrow R^k$ is  Lebesgue measurable, and we call a control
{\it admissible}   if it is ordinary and $u(t)\in U$ for every $t\in[0,t_{f}]$.
Let $L:R^n\times U\rightarrow R$ be an absolutely-integrable
cost function, and let
\begin{equation}
J\ :=\ \int_{0}^{t_{f}}L\big(x(t),u(t)\big)dt
\end{equation}
be its related performance functional. The optimal control problem that we consider is to minimize $J$ over the space of
admissible controls.

The following assumption will be made throughout the paper:
\begin{assumption}
(1).  The function $f(x,u)$ is twice-continuously differentiable in $x\in R^n$ for every $u\in U$;
the functions $f(x,u)$,  $\frac{\partial f}{\partial x}(x,u)$, and
$\frac{\partial^2 f}{\partial x^2}(x,u)$ are locally-Lipschitz continuous in $(x,u)\in R^n\times U$; and
there exists $K>0$ such that, for every $x\in R^n$ and for every $u\in U$,
$||f(x,u)||\leq K(||x||+1)$.
(2). The function $L(x,u)$ is continuously differentiable in $x\in R^n$ for every
$u\in U$; and the functions $L(x,u)$ and   $\frac{\partial L}{\partial x}(x,u)$ are locally-Lipschitz continuous in
$(x,u)\in R^n\times U$.
\end{assumption}
Observe that part (1) of the assumption guarantees the existence of a unique absolutely-continuous solution of
Equation (1) in its integral form for every admissible control and $x_{0}\in R^n$, while part (2) implies that the Lebesgue integral in Equation (2) is well defined.

This paper has been motivated by two kinds of problems: one concerns switched-mode hybrid systems, and the other concerns optimal balancing of the energy required for transmission and motion in mobile robotic networks.
We propose an algorithm defined in the setting of relaxed controls, and analyze its convergence using Polak's framework
of optimality functions. The theoretical results are first derived in the   abstract framework of Eqs. (1) and (2),
and then applied to problems in the two aforementioned areas of interest.

 A standard basic requirement of algorithms for nonlinear-programming (finite-dimensional optimization)
problems, and especially gradient-descent
methods,
is that every accumulation point of an iterate-sequence they compute must satisfy an optimality condition
for a  local minimum, such as the  Kuhn-Tucker condition (e.g., \cite{Polak97}).
Thus, if a bounded sequence of iteration points is computed then it has at least one accumulation point which, therefore,
 must satisfy that condition.
However, in infinite-dimensional optimization this requirement can be  vacuous since bounded sequences
need not have accumulation points.  This issue is not only theoretical but also has practical
implications. An infinite-dimensional problem may not have a solution or  a local minimum on a closed
and bounded set, and even if a local minimum exists,
it is not guaranteed that it can be approximated by the solution points of the problem's finite-dimensional discretizations
at any level of precision.  To get around these issues,
E. Polak developed a comprehensive framework for design and analysis of algorithms for infinite-dimensional
optimization that gives the user considerable discretion in choosing the discretization levels in an adaptive
fashion \cite{Polak97}. A survey of the  framework will be carried out in the next section,
and we mention that recently it has been used in the context of switched-mode  hybrid systems in Refs.
 \cite{Caldwell10,Caldwell11,Caldwell12,Caldwell16,Gonzalez10,Vasudevan12,Egerstedt06,Axelsson08}.

A relaxed control is a mapping $\mu(t)$ from the time interval
$[0,t_f]$ into the space of probability measures on the
set $U$  \cite{McShane67,Warga72}.  Relaxed controls provide a useful framework for optimal control problems for the following
reasons. First, the space  of relaxed controls is convex even though
the input-constraint set $U$ may not be convex. Second, the space of relaxed controls is compact in a suitable topology,
namely the weak-star topology, and hence it contains solutions to optimal control problems lacking solutions in
a functional space of admissible controls like $L^1(U)$ (a more detailed survey of these points is provided in Section 2).
However, a relaxed control is a more abstract object than an ordinary control,
 which can make it problematic for an algorithm to handle an iterative sequence of relaxed controls.  This point will be discussed
 in the sequel.

This paper  combines  the frameworks of optimality functions and
relaxed controls to define a new algorithm for the optimal control problem.  Its main innovation is
in the choice of the search direction from a given relaxed control,  which is based on a pointwise
 minimizer of the  Hamiltonian (defined below) at each time $t\in[0,t_{f}]$.\footnote{The computation of the minimizer of the Hamiltonian does not require
the solving of a two-point boundary value problem, but rather is based on {\it sequential} integrations
of the state equation forwards and the costate (adjoint) equation backwards.
The structure of the problem, and especially the absence of constraints on the final state, make this possible
under Assumption 1.1. Final-state constraints can be handled by the application of penalty functions thereby transforming the problem into
one without final-state constraints, as a forthcoming example will
illustrate.}
Its step size is determined by the Armijo procedure \cite{Armijo66,Polak97}.
To our knowledge it is the first such general algorithm
which is defined entirely in the space of relaxed controls without projecting the results of each iteration into
the space of ordinary controls. The following results will be established.
\begin{itemize}
\item
The aforementioned search direction yields descent of the cost functional (2) even though it is not defined in
explicit gradient terms. It is a form of gradient projection. For a  class of problems  (including autonomous switched-mode systems) its computation is  straightforward and simpler than that of direct gradients.
\item
The algorithm is stable in the sense that it yields descent in the performance integral regardless of
the initial guess, Furthermore,
the simulation results presented in Section 4 indicate its rapid descent  at the early stages of its runs. By this we do not
claim a  high rate of asymptotic convergence since it is a first-order algorithm,
but rather that most of its descent is obtained in few early iterations requiring meagre
computing (CPU) times.
\item
The combination of concepts and techniques from the settings of optimality functions and  relaxed controls yields an analysis framework
that is based on simple arguments. This point will become evident from the simplicity and brevity of the forthcoming proofs.
\end{itemize}

Whereas the area of numerical techniques for optimal control has had a long history (e.g.,  \cite{Polak97} and references therein),
recently there has been a considerable interest in optimal control problems defined on switched-mode hybrid systems.
In this context the problem was formulated   in \cite{Branicky98, Piccoli98}, and variants of the Maximum Principle were derived for it
in \cite{Piccoli98,Sussmann99,Shaikh02,Garavello05,Bengea05,Shaikh07,Taringoo13,Taringoo13a}.
New and emerging algorithmic approaches
include first- and second-order techniques \cite{Xu02,Xu02a,Shaikh02,Egerstedt06,Caldwell11,Vasudevan12},
zoning algorithms based on the geometric properties of
the underlying systems \cite{Shaikh03,Caines05,Shaikh05,Shaikh07},   projection-based algorithms
\cite{Caldwell10,Caldwell12,Caldwell12a,Caldwell12b,Vasudevan12,Caldwell16}, methods based on dynamic programming and convex optimization
\cite{Hellund99},   and algorithms  based
on needle variations \cite{Attia05,Egerstedt06,Axelsson08,Gonzalez10,Taringoo13a}. Concerning the relaxed hybrid problem, Ref. \cite{Ge75} developed generalized-linear programming techniques and convex-programming algorithms,  \cite{Bengea05} derived optimality conditions (both necessary and sufficient) for a class of hybrid optimal control problems, and \cite{Meyer12} applied to them the MATLAB
   \texttt{fmincon} nonlinear programming solver.
Comprehensive recent surveys can be found in  \cite{Zhu11,Lin14}.

The forthcoming algorithm will be presented and analyzed in the abstract problem formulation of  Eqs. (1) and (2) and their extensions
to the relaxed-control setting. However, for implementation, we restrict the class of problems in the following two ways:
1). The pointwise minimizer of the Hamiltonian can be computed or adequately estimated by a simple formula.
2). For every $x\in R^n$,  the dynamic response function $f(x,u)$ in Eq. (1) is affine  in $u$, and the cost function $L(x,u)$ in Eq. (2) is convex in $u$. Many problems of theoretical and practical interest in optimal control satisfy these restrictions. These include
problems defined on autonomous switched-mode systems and other hybrid systems, which will be shown to admit
efficient implementations of the algorithm.  Other kinds of hybrid
systems are not yet included, and these will be mentioned in the sequel as a subject of current research.

The rest of the paper is organized as follows.
Section 2 presents brief surveys of the frameworks of relaxed controls  and optimality functions.
Section 3 describes the algorithm and derives related theoretical results, while Section 4 presents simulation
results.  Finally, Section 5 concludes the paper and points out directions for future research.

${\bf Notation:}$ The term {\it control} refers to the function $\{u:[0,t_{f}]\rightarrow U\}$ and is denoted by the
boldface symbol ${\bf u}$ to distinguish it from a point in the set $U$ which is denoted by the lower-case $u$ or $u(t)$. Similarly, boldface notation
refers to a function of $t\in[0,t_{f}]$ as in ${\bf x}:=\{x:[0,t_{f}]\rightarrow R^n\}$ for the state trajectory,
${\bf p}:=\{p:[0,t_{f}]\rightarrow R^n\}$ for the costate (adjoint) trajectory,
$\boldsymbol{\mu}$ for a relaxed control defined as a function from $t\in[0,t_{f}]$ into the space of probability measures on $U$, etc.

\section{Review of Established Results}
This section recounts the basic framework of relaxed controls and some fundamental notions of algorithms' convergence in infinite-dimensional spaces.

\subsection{Relaxed Control}
Comprehensive presentations of the theory of relaxed controls and their role in optimal control can be found in
 \cite{McShane67,Warga72,Young69,Gamkrelidze78,Berkovitz13}; also see \cite{Lou09} for a recent survey.
 In the following paragraphs we summarize its main points that are relevant to our discussion.
 Let $M$ denote the space of Borel  probability measures on the set $U$, and denote by
 $\mu$ a particular point in $M$.
A relaxed control  associated with the system (1) is a mapping $\mu:[0,t_{f}]\rightarrow M$ which is measurable
in the following sense: For every
continuous function $\phi:U\rightarrow R$, the function $\int_{U}\phi(u)d\mu(t)$ is Lebesgue measurable
in $t$.
We denote
the space of relaxed controls by $\bf{M}$, and   in accordance with previous notation we denote a relaxed control
$\{\mu(t):t\in[0,t_{f}]\}$ by $\boldsymbol{\mu}$.

Recall that an ordinary  control ${\bf u}$ is admissible if the function $u:[0,t_{f}]\rightarrow U$ is (Lebesgue) measurable.
Note that the space of ordinary controls  is embedded in the space of relaxed controls by associating with $u(t)$  the
Dirac
probability measure at $u(t)$ $\forall t\in[0,t_{f}]$. In this case we say,
  with a slight abuse of notation,
  that $\boldsymbol{\mu}$ is an  ordinary control, and indicate this by the notation $\boldsymbol{\mu}\thicksim{\bf u}$.
Furthermore, the space of ordinary controls is dense in the space of relaxed controls in the weak-star topology, namely
 in the following sense: For every relaxed control $\boldsymbol{\mu}$ there exists a sequence
$\{{\bf u}_{k}\}_{k=1}^{\infty}$
of ordinary controls such that, for every function $\psi\in L^{1}\big([0,t_{f}];C(U)\big)$,\footnote{$L^{1}\big([0,t_{f}];C(U)\big)$
is the space of functions $\psi:[0,t_{f}]\times U\rightarrow R$ that are measurable and absolutely integrable
in $t$ on $[0,t_{f}]$ for every $u\in U$, and
continuous on $U$ for every $t\in[0,t_{f}]$.}
\[
\lim_{k\rightarrow\infty}\int_{0}^{t_{f}}\psi(t,u_{k}(t))dt\ =\ \int_{0}^{t_{f}}\int_{U}\psi(t,u)d\mu(t)dt.
\]
 Furthermore, the space of relaxed controls is compact in the weak-star topology.

An extension of the system defined by Equations (1) and (2) to the setting of relaxed controls is provided by the state equation
\begin{equation}
\dot{x}(t)=\int_{U}f\big(x(t),u\big)d\mu(t)
\end{equation}
with the same boundary condition $x_{0}=x(0)$ as for (1),
and the related cost functional defined as
\begin{equation}
J(\boldsymbol{\mu})=\int_{0}^{t_{f}}\int_{U}L\big(x(t),u\big)d\mu(t) dt.
\end{equation}
The relaxed optimal control problem is to minimize $J(\boldsymbol{\mu})$ over $\boldsymbol{\mu}\in{\bf M}$.

There are two noteworthy special cases. First, if $\boldsymbol{\mu}\thicksim{\bf u}$,  then Equations (3) and (4) are reduced to Equations (1) and (2), respectively.
Second, in the case where $U:=\{u_{1},\ldots,u_{m}\}$ is a finite set, (3) and (4) have the following respective forms:
$
\dot{x}(t)=\sum_{i=1}^{m}\mu^{i}(t)f(x(t),u_{i})
$
 and $J=\int_{0}^{t_{f}}\sum_{i=1}^{m}\mu^{i}(t)L(x(t),u_{i})dt$, with $\mu^{i}(t)\geq 0$ and $\sum_{i=1}^{m}\mu^{i}(t)=1$
 $\forall t\in[0,t_{f}]$,
and this corresponds to the case of autonomous switched-mode systems.  The space of relaxed controls
generally is convex even though the set $U$ need not be convex.

Essential parts of the theory of optimal control, including the Maximum Principle \cite{Vinter00}, apply to the relaxed-control problem
(see also \cite{Young69,McShane67,Warga72,Gamkrelidze78,Berkovitz13}). Thus,
defining the adjoint (costate) variable $p(t)\in R^n$ by the equation
\begin{equation}
\dot{p}(t)=-\int_{U}\Big(\frac{\partial f}{\partial x}\big(x(t),u\big)^{\top}p(t)+\frac{\partial L}{\partial x}\big(x(t),u\big)^{\top}\Big)d\mu(t)
\end{equation}
with the boundary condition $p(t_{f})=0$,  the Hamiltonian has the form
\begin{align}
H\big(x(t),\mu(t),p(t)\big) \nonumber \\
=\int_{U}\Big(p(t)^{\top}f\big(x(t),u\big)+L\big(x(t),u\big)\Big)d\mu(t),
\end{align}
and the
 Maximum Principle states that if $\boldsymbol{\mu}\in{\bf M}$ is a minimum for the relaxed optimal control problem then
$\mu(t)$ minimizes the Hamiltonian at almost every time-point $t\in[0,t_{f}]$.\footnote{Note that  the integrand
of (6) is the usual Hamiltonian $H(x,u,p):=p^{\top}f(x,u)+L(x,u)$, while the term in the Left-hand Side of (6), namely
$H(x,\mu,p)$, refers to the relaxed Hamiltonian. These two notations are distinguished by their second variable, $u$ vs. $\mu$,
which suffices to  render the usage of the functional term
 $H(x,\cdot,p)$ unambiguous  in the sequel.
}

\subsection{Infinite-dimensional Optimization}
It is a common
practice to characterize convergence of algorithms for infinite-dimensional
optimization problems in terms of {\it optimality functions} \cite{Polak97}.
 Consider the abstract optimization problem of minimizing a function $\phi:\Gamma\rightarrow R$
where $\Gamma$ is a topological space, and consider an optimality condition (necessary or sufficient) associated with this optimization
problem. An optimality function is a function $\theta:\Gamma\rightarrow R^-$ having the property that
$\theta(v)=0$ if and only if $v$  satisfies the
 optimality condition. The optimality-function concept is useful if $|\theta(v)|$ is
 a meaningful heuristic measure of the extent to which $v\in\Gamma$ fails to satisfy the optimality condition. For example,
 if $\phi$ is a continuously Frechet-differentiable functional defined on a Hilbert space $H$, then an optimality condition
 is $\frac{d\phi}{dv}(v)=0$ ($\frac{d\phi}{dv}$ meaning the Frechet derivative), and a meaningful associated optimality function
 is $\theta(v):=-||\frac{d\phi}{dv}(v)||$, where the indicated norm is in $H$.

Reference \cite{Polak97} devised a framework for analysis of algorithms in this abstract
setting,  where convergence  is defined as follows: If an
algorithm computes a sequence $v_{k}$, $k=1,2,\ldots$, of points in $\Gamma$ then
\begin{equation}
\lim_{k\rightarrow\infty}\theta(v_{k})=0.
 \end{equation}
 In this abstract setting such  a sequence $\{v_{k}\}_{k=1}^{\infty}$ need not have an accumulation point
 even if it is a bounded sequence in a metric space (unless it is isomorphic to a Euclidean space),  and therefore a characterization of
 an algorithm's convergence in terms of such accumulation points could be vacuous. The use of optimality functions
 via Equation  (7) serves to  resolve this conceptual issue. Note that it is a form of weak convergence.

 Consider  an  algorithm that computes from a given $v\in\Gamma$ the next iteration point,
 denoted by $v_{{\rm next}}$, and suppose that its repetitive application
 computes a sequence of iteration points  $\{v_{k}\}_{k=1}^{\infty}\subset\Gamma$,
 where $v_{k+1}=v_{k,{\rm next}}$. The
 algorithm
is said to have the property of {\it sufficient descent} with respect to $\theta(\cdot)$ if the following two conditions hold:
(i) $\phi(v_{{\rm next}})-\phi(v)\leq 0$ for every $v\in\Gamma$, and (ii) for every $\eta>0$ there exists
$\delta>0$ such that for every $v\in\Gamma$, if $\theta(v)<-\eta$ then
$\phi(v_{{\rm next}})-\phi(v)<-\delta$. For such an algorithm,
the following result is a straightforward  corollary of Theorem 1.2.8 in
\cite{Polak97}  and hence its proof is omitted.
\begin{proposition}
Suppose that $|\phi(v)|$ is bounded over $\Gamma$. If an algorithm is of sufficient descent then it is convergent in the sense of (7).\hfill$\Box$
\end{proposition}
Refs.
\cite{Caldwell10,Caldwell11,Caldwell12,Caldwell16,Gonzalez10,Vasudevan12,Egerstedt06,Axelsson08} used the framework of
optimality functions and sufficient descent to define and analyze their respective algorithms for the switched-mode
 optimal control problem,
 where $\Gamma$ is the space of admissible controls ${\bf u}$, and
$\theta({\bf u})$ typically is related to the magnitude of the steepest feasible descent-direction vector. In contrast,
the optimality function defined in this paper is based on  the Hamiltonian rather than  the steepest descent or  any explicit form of
a derivative.

Consider the relaxed control problem defined by
Equations (3) and (4). Given a relaxed control $\boldsymbol{\mu}\in{\bf M}$, let ${\bf x}$  and ${\bf p}$ be the
associated state trajectory and costate trajectory
as defined by  (3) and (5), respectively. We use the following optimality function, $\theta(\boldsymbol{\mu})$:\footnote{In this and later equations we drop the explicit notational dependence of various integrand-terms  on $t$ when no confusion arises.}
\begin{equation}
\theta(\boldsymbol{\mu})\ =\ \min_{\boldsymbol{\nu}\in{\bf M}}\int_{0}^{t_{f}}\big(H(x,\nu,p)-H(x,\mu,p)\big)dt,
\end{equation}
where the Hamiltonians in the integrand of (8) were defined in (6).  Recall that $x$ and $p$ in (8) are associated with $\mu$ and hence are independent of $\nu$; therefore, the compactness of the space of relaxed controls
implies that the minimum (not only {\rm inf})
in (8) exists.
Let $\boldsymbol{\mu}^{\star}\in{\bf M}$ denote an {\rm argmin}, then (8) becomes
\begin{equation}
\theta(\boldsymbol{\mu})=\int_{0}^{t_{f}}\big(H(x,\mu^{\star},p)-H(x,\mu,p)\big)dt.
\end{equation}
We observe that this optimality function satisfies the aforementioned properties  with respect to the Maximum Principle: Obviously
 $\theta(\boldsymbol{\mu})\leq 0$
 for every $\boldsymbol{\mu}\in{\bf M}$;
 $\theta(\boldsymbol{\mu})=0$ if and only if $\boldsymbol{\mu}$
minimizes the Hamiltonian at almost every $t\in[0,t_{f}]$ and hence  satisfies the
Maximum Principle; and $|\theta(\boldsymbol{\mu})|$ arguably indicates the extent to which the Maximum Principle is not satisfied at $\boldsymbol{\mu}$.

\section{Hamiltonian-based Algorithm}

The analysis in this section is carried out under Assumption 1.1.

Given a relaxed control $\boldsymbol{\mu}\in{\bf M}$, let ${\bf x}$ and ${\bf p}$ denote the
related state trajectory and costate trajectory defined by (3) and (5), respectively. For   these
state and costate, and for every $t\in[0,t_{f}]$,
consider the Hamiltonian $H(x(t),\cdot,p(t)\big)$, defined by (6), as a function of its second variable.
Fix another relaxed control $\boldsymbol{\nu}\in{\bf M}$. Now for every $\lambda\in[0,1]$, $\lambda\boldsymbol{\nu}+(1-\lambda)\boldsymbol{\mu}$
is also a  relaxed control, and we denote it by
$\boldsymbol{\mu}_{\lambda}$. Furthermore, let $\{x_{\lambda}(t):\ \ t\in[0,t_{f}]\}$ denote the state trajectory
associated with $\boldsymbol{\mu}_{\lambda}$ as defined by (3), namely,
\begin{align}
\dot{x}_{\lambda}(t\ )=\ \lambda\int_{U}f\big(x_{\lambda}(t),u\big)d\nu(t)\nonumber \\
+(1-\lambda)\int_{U}f\big(x_{\lambda}(t),u\big)d\mu(t),
\end{align}
and define  $\tilde{J}(\lambda):=J(\boldsymbol{\mu}_{\lambda})$, where $J$ is defined by (4),
namely
\begin{align}
\tilde{J}(\lambda) & =\int_{0}^{t_{f}}\Big(\lambda
\int_{U}L\big(x_{\lambda}(t),u\big)d\nu(t) \nonumber \\
 & +(1-\lambda)\int_{U}L\big(x_{\lambda}(t),u\big)d\mu(t)\Big)dt.
\end{align}
The algorithm described in this section is based on moving from  $\boldsymbol{\mu}\in{\bf M}$ in the direction
of  $\boldsymbol{\nu}$ by choosing a step size $\lambda\in[0,1]$, and therefore, we next characterize those
$\boldsymbol{\nu}\in{\bf M}$ that provide a direction of descent.
\begin{proposition}
The one-sided derivative $\frac{d\tilde{J}}{d\lambda^{+}}(0)$ exists and has the following form,
\begin{equation}
\frac{d\tilde{J}}{d\lambda^{+}}(0)\ =\ \int_{0}^{t_{f}}\big(H(x,\nu,p)-H(x,\mu,p)\big)dt,
\end{equation}
where all the terms in the integrand in the Right-Hand Side (RHS) of (12) are functions of time.
\end{proposition}
{\it Proof.}
Consider the Right-Hand Sides (RHS) of Equations (10) and (11) as functions of $x=x_{\lambda}(t)$ and $\lambda\in[0,1]$,
for given relaxed controls $\boldsymbol{\mu}$ and $\boldsymbol{\nu}$. By Assumption 1.1 these functions are
twice-continuously differentiable ($C^2$) in $x$, and by (10) and (11)  they are linear in $\lambda$ and hence
$C^2$ as well. Therefore,
standard variational techniques show  that  the function $\tilde{J}(\lambda)$ is differentiable on
$\lambda\in[0,1]$ and its derivative has the following form,
\begin{align}
\frac{d\tilde{J}}{d\lambda}(\lambda) & =\int_{0}^{t_{f}}\Big(p_{\lambda}(t)^{\top}\int_{U}f\big(x_{\lambda}(t),u\big)\big(d\nu(t)-
d\mu(t)\big) \nonumber \\
& +
\int_{U}L\big(x_{\lambda}(t),u\big)\big(d\nu(t)-d\mu(t)\big)\Big)dt,
\end{align}
where $p_{\lambda}(t)$ is the costate associated with this derivative. Furthermore, by (10), it is  seen that
${\bf p}_{\lambda}$ is given by Equation (5) with $\boldsymbol{\mu}_{\lambda}$ instead of $\boldsymbol{\mu}$.

Next, for $\lambda=0$, we have that $\boldsymbol{\mu}_{0}=\boldsymbol{\mu}$, ${\bf x}_{0}=\bf{x}$, and ${\bf p}_{0}={\bf p}$,
and therefore, with $\lambda=0$ in (13) we obtain,
\begin{align}
\frac{d\tilde{J}}{d\lambda^{+}}(0) & =\int_{0}^{t_{f}}\Big(p(t)^{\top}\int_{U}f\big(x(t),u\big)\big(d\nu(t)-d\mu(t)\big) \nonumber \\
& +
\int_{U}L\big(x(t),u\big)\big(d\nu(t)-d\mu(t)\big)\Big)dt.
\end{align}
By Equation (6) the RHS of (14) is identical to the RHS of (12).
\hfil\hfill $\Box$

Equation (12) implies that $\boldsymbol{\nu}$ is a descent direction from $\boldsymbol{\mu}$ if the RHS of (12) is negative.
This is the case if $\nu(t)$ is a pointwise minimizer of the Hamiltonian over $M$ at each time $t\in[0,t_{f}]$, unless
$H(x(t),\nu(t),p(t))=H(x(t),\mu(t),p(t))$ for almost every $t\in[0,t_{f}]$. In this case,
let us use the notation
$\boldsymbol{\mu}=\boldsymbol{\nu}$ for a pointwise minimizer of the Hamiltonian.
The following result implies that the pointwise search for
such a minimizer can be confined to $U$ and need not be extended to $M$, the space of Borel
probability measures on $U$.
\begin{proposition}
Fix $x\in R^n$ and $p\in R^n$.  Let $u^{\star}\in{\rm argmin}\{(H(x,u,p):u\in U\}$. Then, for every probability measure
$\nu\in M$,
\begin{equation}
H(x,u^{\star},p)\leq H(x,\nu,p).
\end{equation}
\end{proposition}
Note that the Left-hand side (LHS) of (15) is the usual Hamiltonian while its RHS is the relaxed Hamiltonian defined by (6).

{\it Proof.}
By (6) and the fact that $u^{\star}$ minimizes the ordinary Hamiltonian over $u\in U$, we have, for
every $\nu\in M$,
\begin{align}
& H(x,\nu,p) =\int_{U}H(x,u,p)d\nu \nonumber \\
& \geq\int_{U}H(x,u^{\star},p)d\nu=H(x,u^{\star},p).
\end{align}
 \hfill$\Box$

 We point out that the  point $u^{\star}$ in the statement of Proposition 3.2 is the pointwise minimizer of the Hamiltonian
over $U$, for given $x\in R^n$ and $p\in R^n$. Such a minimizer exists and the minimum is finite since,
by assumption the set $U$ is
compact, and by Assumption 1.1 the function $H(x,u,p)$ is continuous in
$u\in U$.

Let $\boldsymbol{\mu}$ be a relaxed control, and let ${\bf x}$ and ${\bf p}$ be the associated state trajectory and costate trajectory
as defined by Equations (3) and (5), respectively.  For every $t\in[0,t_{f}]$, let $u^{\star}(t)\in{\rm argmin}\{H(x(t),u,p(t)):u\in U\}$.
It does not mean that the function $\{u^{\star}(t),\ t\in[0,t_{f}]\}$, is an admissible control
since it might not be Lebesgue measurable. On the other hand, we have seen that there exists a relaxed control $\boldsymbol{\mu}^{\star}$ that minimizes the RHS of (8) over $\boldsymbol{\nu}\in{\bf M}$, and hence
$\mu^{\star}(t)$ minimizes $H\big(x(t),\nu,p(t)\big)$ over $\nu\in M$ (the space of Borel probability measures on $U$)
for almost every $t\in[0,t_{f}]$.

Ideally we would like to choose such $\boldsymbol{\mu}^{\star}$ as the descent direction of the algorithm from
$\boldsymbol{\mu}$, but its computation may be fraught with difficulties for the following two reasons: (i) for a given
$t$, $\mu^{\star}(t)$ may
 not be a Dirac measure at a point in $U$, and
(ii) The pointwise minimizer $\mu^{\star}(t)$ has to be computed for every $t$ in the infinite set $[0,t_{f}]$.
Therefore we choose as descent direction a relaxed control $\boldsymbol{\nu}\in{\bf M}$ having the
following two properties: (i)
$\boldsymbol{\nu}\thicksim{\bf v}$ where
${\bf v}$ is a piecewise-constant ordinary control, and (ii)
$\boldsymbol{\nu}\in{\bf M}$
 ``almost'' minimizes the Hamiltonian in the
following sense: For a given a constant $\eta\in(0,1)$ which we fix throughout the algorithm
(below),
\begin{equation}
\int_{0}^{t_{f}}\big(H(x,\nu,p)-H(x,\mu,p)\big)dt\leq\eta\theta(\boldsymbol{\mu});
\end{equation}
$x,\ \nu,\ p,$ and $\mu$ are all functions of time.
We label such $\boldsymbol{\nu}$ an {\it $\eta$-minimizer of the Hamiltonian}.
It will be seen that it is always possible to choose a relaxed control $\boldsymbol{\nu}$ with these two properties
as long as $\theta(\boldsymbol{\mu})<0$.
With this direction of descent, the algorithm uses the Armijo step size
\cite{Armijo66,Polak97}. It has the following form.

Given constants $\alpha\in(0,1)$, $\beta\in(0,1)$, and $\eta\in(0,1)$.
\begin{algo}
Given $\boldsymbol{\mu}\in{\bf M}$,
compute $\boldsymbol{\mu}_{{\rm next}}\in{\bf M}$ by the following steps.

{\it Step 0:} If $\theta(\boldsymbol{\mu})=0$, set $\boldsymbol{\mu}_{{\rm next}}=\boldsymbol{\mu}$, then exit.\\
{\it Step 1:} Compute the state and costate trajectories, ${\bf x}$ and ${\bf p}$,
associated with $\boldsymbol{\mu}$,
by using Equations (3) and (5), respectively.\\
{\it Step 2:} Compute a relaxed control $\boldsymbol{\nu}\in{\bf M}$ which is an $\eta$-minimizer of the Hamiltonian, namely it
satisfies Equation (17).\\
 {\it Step 3:} Compute the integer $\ell_{\mu}$ defined as follows,
 \begin{align}
 & \ell_{\mu}=\min\big\{\ell=0,1,\ldots,: \nonumber \\
 & J(\boldsymbol{\mu}+\beta^{\ell}(\boldsymbol{\nu}-\boldsymbol{\mu}))-J(\boldsymbol{\mu})\leq\alpha
 \beta^{\ell}\eta\theta(\boldsymbol{\mu})\big\}.
 \end{align}
 Define  $\lambda_{\mu}:=\beta^{\ell_{\mu}}$.\\
 {\it Step 4:} Set
 \begin{equation}
  \boldsymbol{\mu}_{{\rm next}}=\boldsymbol{\mu}+\lambda_{\mu}(\boldsymbol{\nu}-\boldsymbol{\mu}).
 \end{equation}
 \end{algo}
 \vspace{.05in}

 A few remarks are due.

 1).  The algorithm is meant to be run iteratively and compute a sequence $\{\boldsymbol{\mu}_{k}\}_{k\geq 1}$
such that $\boldsymbol{\mu}_{k+1}=\boldsymbol{\mu}_{k,{\rm next}}$, as long as it does not exit in   Step 0.

2). The algorithm does not attempt to solve a two-point boundary value problem. In Step 1 it first integrates the
differential equation (3) forward from the initial condition $x_{0}$, and then the differential equation
(5) backwards from the specified terminal condition $p(t_{f})=0$.   We do not specify the particular numerical
integration technique that should be used, but say more about it in the sequel.

3). We do not specify the choice of  $\boldsymbol{\nu}$ in Step 2 but rather leave it to the user's
discretion. However,  we point out  that such an $\eta$-minimizer of the Hamiltonian always exists in the form of
$\boldsymbol{\nu}\thicksim{\bf v}$ for a piecewise-constant ordinary control
${\bf v}$, unless $\theta(\boldsymbol{\mu})=0$. The reason is that
  the space of ordinary controls is dense in the space of relaxed controls in the weak star topology on
 $L^{1}\big([0,t_{f}];C(U)\big)$, and the space of piecewise-constant ordinary controls is dense in the space of ordinary controls in the $L^1$ norm and hence in the weak-star topology.

 4). In Step 4,  $\mu_{{\rm next}}(t)$ is a
 convex combination of $\boldsymbol{\mu}$ and $\boldsymbol{\nu}$
  in the sense of measures, meaning that
 for every continuous function $g:U\rightarrow R$,
 \[
 \int_{U}g(u)d\mu_{{\rm next}}(t)=\int_{U}g(u)\big((1-\lambda_{\mu})d\mu(t)+\lambda_{\mu}d\nu(t)\big).
 \]
 In the event that  $\boldsymbol{\mu}\thicksim{\bf u}$ and $\boldsymbol{\nu}\thicksim{\bf v}$, this means that
 \[
 \int_{U}g(u)d\mu_{{\rm next}}(t)=(1-\lambda_{\mu})g(u(t))+\lambda_{\mu}g(v(t)),
 \]
 which does not necessarily imply that
 $\boldsymbol{\mu}_{{\rm next}}\thicksim{\bf u}+\lambda_{\mu}({\bf v}-{\bf u})$.

 5). The Armijo step size, computed in Step 3, is commonly used in nonlinear programming as well as in infinite-dimensional optimization.
 Reference \cite{Polak97} contains analyses of several algorithms using it and practical guidelines
 for its implementation.

 6). It will be proven that the integer $\ell_{\mu}$ defined in Step 3 is finite as long as $\theta(\boldsymbol{\mu})<0$, hence the
 algorithm cannot jam at a point (relaxed control) that does not satisfy the Maximum Principle, namely the condition
 $\theta(\boldsymbol{\mu})=0$.

We point out that the idea of a descent direction comprised of a pointwise minimizer of the Hamiltonian has its origin in \cite{Mayne75}.
The algorithm proposed in that reference is defined in the setting of ordinary
 controls, its descent direction is all the way to
a minimizer of the Hamiltonian
at a subset of the time-horizon $[0,t_{f}]$, and its main convergence result is stated in terms
of  accumulation points of  computed iterate-sequences.  The algorithm in this paper is quite different in that
it is defined in the setting of relaxed controls, it moves part of the way towards a minimizer of the Hamiltonian throughout the
entire
time horizon, and its analysis is carried out in the context of optimality functions.

 We next establish convergence of Algorithm 3.3.

 \begin{proposition}
 Suppose that Assumption 1.1 is satisfied. Let $\{\boldsymbol{\mu}_{k}\}_{k=1}^{\infty}$ be a sequence of
 relaxed controls computed by Algorithm 3.3, such that for every $k=1,2,\ldots$,  $\boldsymbol{\mu}_{k+1}=\boldsymbol{\mu}_{k,next}$. Then,
 \begin{equation}
 \lim_{k\rightarrow\infty}\theta(\boldsymbol{\mu}_{k})=0.
 \end{equation}
 \end{proposition}
 The proof is based on the following lemma, variants of which have been proved
 in \cite{Polak97} (e.g., Theorem 1.3.7). We supply the proof in order to complete the presentation.
\begin{lemma}
Let $g(\lambda):R\rightarrow R$ be a twice-continuously differentiable $(C^{2}$) function. Suppose that
$g^{\prime}(0)\leq 0$, and there exists $K>0$
 such that
$|g^{\prime\prime}(\lambda)|\leq K$ for every $\lambda\in R$.
Fix $\alpha\in(0,1)$, and define $\gamma:=2(1-\alpha)/K$. Then for  every positive $\lambda\leq\gamma|g^{\prime}(0)|$,
\begin{equation}
g(\lambda)-g(0)\ \leq\ \alpha\lambda g^{\prime}(0).
\end{equation}
\end{lemma}
{\it Proof.}
Recall (see \cite{Polak97}, Eq. (18b), p.660) the following   exact second-order expansion    of $C^2$ functions,
\begin{equation}
g(\lambda)\ =\ g(0)+\lambda g^{\prime}(0)+\lambda^2\int_{0}^{1}(1-s)g^{\prime\prime}(s\lambda)ds.
\end{equation}
Using this and the assumption that $|g^{\prime\prime}(\cdot)|\leq K$, we obtain that
\begin{align}
& g(\lambda)-g(0)-\alpha\lambda g^{\prime}(0)\nonumber \\
=\ & \ (1-\alpha)\lambda g^{\prime}(0)+\lambda^2\int_{0}^{1}(1-s)g^{\prime\prime}(s\lambda)ds   \nonumber \\
\leq\ & (1-\alpha)\lambda g^{\prime}(0)+\lambda^2K/2\nonumber \\
=\ & \ \lambda\big((1-\alpha)g^{\prime}(0)+\lambda K/2\big).
\end{align}
For every positive $\lambda\leq\gamma|g^{\prime}(0)|$, $(1-\alpha)g^{\prime}(0)+\lambda K/2\leq 0$, and hence, and by (23), Equation (21) follows.
\hfill $\Box$

{\it Proof of Proposition 3.4.}\ \
Let $\boldsymbol{\mu}\in{\bf M}$ and $\boldsymbol{\nu}\in{\bf M}$ be any two relaxed controls, and
for $\lambda\in[0,1]$, consider $\tilde{J}(\lambda)$ as defined by Equation (11). We next show that $\tilde{J}$ is a
twice-continuously differentiable function of $\lambda$, and there exists $K>0$ such that, for all
$\boldsymbol{\mu}\in{\bf M}$,
$\boldsymbol{\nu}\in{\bf M}$, and $\lambda\in[0,1]$,
\begin{equation}
\Big|\frac{d^{2}\tilde{J}}{d\lambda^{2}}(\lambda)\Big|\leq K.
\end{equation}
This follows from variational arguments
developed and summarized in \cite{Polak97} as follows. First,
consider two ordinary controls, ${\bf u}$ and ${\bf w}$, and define
${\bf u}_{\lambda}:=\lambda{\bf w}+(1-\lambda){\bf u}$ for every $\lambda\in[0,1]$. Denote by ${\bf x}_{\lambda}$
the state trajectory associated with ${\bf u}_{\lambda}$ as defined by (1), and let $\tilde{J}(\lambda):=J({\bf u}_{\lambda})$ be the cost functional,
defined by (2), as a function of $\lambda$.  By Assumption 1.1,
an application of Corollary 5.6.9 and Proposition 5.6.10
in
\cite{Polak97} to variations in $\lambda$ yields that $\tilde{J}(\lambda)$ is continuously differentiable in $\lambda$,
and  the term
$|\frac{d\tilde{J}}{d\lambda}(\lambda)|$ is  bounded from above over all ordinary
admissible controls ${\bf u}$, ${\bf w}$, and $\lambda\in[0,1]$. A second application of these   arguments to the derivative $\frac{d\tilde{J}}{d\lambda}$, supported by
the $C^2$ assumption (Assumption 1.1),  yields
that $\tilde{J}(\lambda)$ is twice continuously differentiable and its second derivative also is  bounded from above over all
ordinary admissible controls ${\bf u}$, ${\bf w}$, and $\lambda\in[0,1]$.

By the Lebesgue Dominated Convergence Theorem, the same result
holds true when ${\bf u}$ and ${\bf w}$ are replaced by two respective
relaxed controls, $\boldsymbol{\mu}$ and $\boldsymbol{\nu}$, and   the state equation and cost function are defined
by Equations (3) and (4), respectively. This shows that  $\tilde{J}$ is $C^2$ in $\lambda$, and there exists $K>0$,
independent of
$\boldsymbol{\mu}\in{\bf M}$,
$\boldsymbol{\nu}\in{\bf M}$, and $\lambda\in[0,1]$, such that Equation (24) is satisfied.

Let us apply this result to $\boldsymbol{\mu}$ and $\boldsymbol{\nu}$ where  $\boldsymbol{\mu}\in{\bf M}$
is a given relaxed control and
$\boldsymbol{\nu}$ is  defined in Step 2 of Algorithm 3.3.
Recall the constant $\alpha\in(0,1)$ that is used by the algorithm, and define
$\gamma:=2(1-\alpha)/K.$ By Lemma 3.5, for every positive $\lambda<\gamma|\frac{d\tilde{J}}{d\lambda^{+}}(0)|,$
\begin{equation}
\tilde{J}(\lambda)-\tilde{J}(0)\leq\alpha\lambda\frac{d\tilde{J}}{d\lambda^{+}}(0).
\end{equation}
Recall that $\tilde{J}(\lambda):=J(\boldsymbol{\mu}_{\lambda})$, and therefore, according to  the notation in the RHS of (18),
$J(\boldsymbol{\mu}+\beta^\ell(\boldsymbol{\mu}^{\star}-\boldsymbol{\mu}))-J(\boldsymbol{\mu})=\tilde{J}(\beta^\ell)-\tilde{J}(0)$;
consequently, if $\beta^\ell<\gamma|\frac{d\tilde{J}}{d\lambda^{+}}(0)|$, then (25) is satisfied with $\lambda=\beta^\ell$, namely,
\begin{equation}
J(\boldsymbol{\mu}+\beta^\ell(\boldsymbol{\mu}^{\star}-\boldsymbol{\mu}))-
J(\boldsymbol{\mu})\leq\alpha\beta^\ell\frac{d\tilde{J}}{d\lambda^{+}}(0).
\end{equation}
By Proposition 3.1 (Equation (12)) and  the choice of $\boldsymbol{\nu}$ in Step 2 (Equation (17)),
 the RHS of (26) implies that
\begin{equation}
J(\boldsymbol{\mu}+\beta^\ell(\boldsymbol{\mu}^{\star}-\boldsymbol{\mu}))- J(\boldsymbol{\mu})\leq\alpha\beta^\ell\eta\theta(\boldsymbol{\mu})
\end{equation}
as long as $\beta^\ell<\gamma|\frac{d\tilde{J}}{d\lambda^{+}}(0)|$.
But (12), (17), and the fact that  $\frac{d\tilde{J}}{d\lambda^{+}}(0)\leq 0$ imply that
$|\frac{d\tilde{J}}{d\lambda^{+}}(0)|\geq\eta|\theta(\boldsymbol{\mu})|$ and hence
$\gamma|\frac{d\tilde{J}}{d\lambda^{+}}(0)|\geq\eta\gamma|\theta(\boldsymbol{\mu})|$;
consequently (27) is satisfied as long as  $\beta^\ell<\eta\gamma|\theta(\boldsymbol{\mu})|$. Therefore, and by Equation (18), the step size $\lambda_{\mu}$ defined in Step
3 satisfies the inequality
\begin{equation}
\lambda_{\mu}:=\beta^{\ell_{\mu}}\geq\beta\eta\gamma|\theta(\boldsymbol{\mu})|.
\end{equation}
Next, Equations (18) and (19) imply that
\begin{equation}
J(\boldsymbol{\mu_{{\rm next}}})-J(\boldsymbol{\mu})\leq\alpha\lambda_{\mu}\eta\theta(\boldsymbol{\mu}),
\end{equation}
and hence, and by (28) and the fact that $\theta(\boldsymbol{\mu})\leq 0$, we have that
\begin{equation}
J(\boldsymbol{\mu_{{\rm next}}})-J(\boldsymbol{\mu})\leq-\alpha\beta\gamma\eta^2\theta(\boldsymbol{\mu})^2.
\end{equation}
But $\gamma$ is independent of $\boldsymbol{\mu}$ or
$\boldsymbol{\nu}$, and hence Equation (30) implies that the algorithm is of a sufficient  descent.

Finally, the set $U$ is compact by assumption, and therefore standard applications of the Bellman-Gronwall inequality and Equation (4)  yield that $|J(\boldsymbol{\mu})|$  is upper-bounded over all
$\boldsymbol{\mu}\in{\bf M}$. Consequently Proposition 2.1 implies the validity of Equation (20). \hfill $\Box$

 {\it A note on implementation.} Implementations of Algorithm 3.3 generally require numerical integration methods
for computing (or approximating) ${\bf x}$ and ${\bf p}$ in Step 1, and a computation of  $\boldsymbol{\nu}$ in Step 2 that
is based on the pointwise
minimization of the Hamiltonian at a finite number of points in the time-horizon $[0,t_{f}]$. Both require finite grids on the time horizon,
which may be different and can vary from one iteration to the next. The choice of the grid sizes generally comprises
a balance between precision and computing times. A rule-of-thumb proposed in \cite{Polak97} is to adjust the grid adaptively
by tightening it whenever it is sensed that a local optimum is approached.
This adaptive-precision technique underscores Polak's algorithmic framework of {\it consistent approximation} for
infinite-dimensional optimization while guaranteeing convergence in the sense of Eq. (20).

In this paper we are not concerned with the formal rules for adjusting the grids (and hence precision levels). Instead,
we  run the algorithm several times per problem, each with a fixed grid, to see how far it can minimize the cost functional.
The goal of this experiment is to test the tradeoff between precision and
computing times. The results, presented in the next section,
indicate rapid descent from the initial guess regardless of how far it
is from the optimum.  In fact, the key argument in the proof of convergence is the {\it sufficient  descent}
property of the algorithm, captured in Equation (30), which implies large descents until an iteration-sequence
approaches a local minimum. This suggests that the main utility of the algorithm is not in the asymptotic convergence
close to a minimum (where higher-order methods can be advantageous), but rather in its approach to such points.
For a detailed discussion of this point and some comparative results please see Section 4.

Another point related to implementation concerns the algorithm's having to compute  relaxed controls,
which are more complicated objects than ordinary controls. To address this concern we next examine a class of systems where
the state equation $f(x,u)$ is affine in $u$  and the cost function $L(x,u)$ is convex in $u$.

Consider the case where
\begin{equation}
f(x,u)=\phi_{f}(x)+\Psi_{f}(x)u,
\end{equation}
where the functions $\phi_{f}:R^n\rightarrow R^n$ and  $\Psi_{f}:R^n\rightarrow R^{n\times k}$ (the latter being
the space of $n\times k$ matrices)
satisfy Assumption 1.1. By the linearity of the integration operator
and the fact that $\mu(t)$ is a probability measure, Equation (3) assumes the form
\begin{equation}
\dot{x}(t)=f\Big(x(t),\int_{U}ud\mu(t)\Big),
\end{equation}
meaning that, for the purpose of computing the state trajectory,
 the convexification of the vector field inherent in (3) yields the same results as
a convexification of the control. Defining $\bar{u}(t):=\int_{U}ud\mu(t)$, the state equation becomes $\dot{x}=f(x,\bar{u})$,
and we can view $\bar{{\bf u}}$ as a control function from $[0,t_{f}]$ into $conv(U)$. Likewise, if
$L(x,u)= \phi_{L}(x)+\psi_{L}(x)u$ with $\phi_{L}:R^n\rightarrow R$ and $\psi_{L}:R^n\rightarrow R^{1\times k}$,
 then (4) becomes
 \begin{equation}
 J(\boldsymbol{\mu})=\int_{0}^{t_{f}}L(x,\bar{u})dt.
 \end{equation}
 In this case the relaxed optimal control problem is cast as an ordinary  optimal
 control problem with the input constraints $\bar{u}(t)\in conv(U)$. Moreover, if Algorithm 3.3 starts at an ordinary
 control on $conv(U)$ then it could compute only
 such ordinary controls. The reason is that $\boldsymbol{\nu}$ in Step
  2 can always be an ordinary control (as earlier said), and  the convexification of measures in Eq. (19)
 can be carried out by  the convexification of ordinary controls. Therefore, if in Eq. (19), $\boldsymbol{\mu}\thicksim{\bf u}$
 and $\boldsymbol{\nu}\thicksim{\bf v}$ for ordinary admissible controls (on $conv(U)$), then (19) yields
 \begin{equation}
 \boldsymbol{\mu}_{{\rm next}}\thicksim{\bf u}+\lambda_{\mu}({\bf v}-{\bf u}),
 \end{equation}
 which is an ordinary admissible control on $conv(U)$. This is the case of autonomous switched-mode systems, as will be
 demonstrated
 in Section 4.

 Consider next the case where $f(x,u)$ is affine in $u$ as in Equation (31), and $L(x,u)$ is convex in $u$ for every $x\in R^n$.
Then Equation (32) is true but (33) and hence (34) are not true. Therefore, supposing that
$\boldsymbol{\mu}\thicksim{\bf u}$
 and $\boldsymbol{\nu}\thicksim{\bf v}$ for ordinary controls (on $conv(U)$), Equation (19) yields that
$\boldsymbol{\mu}_{{\rm next}}$
is a relaxed control but not an ordinary control on $conv(U)$. However, the convexity of $L(x,u)$ in $u$ in conjunction
with  (32)  imply  that, for every $\lambda\in[0,1]$
\begin{eqnarray}
J\big({\bf u}+\lambda({\bf v}-{\bf u})\big)\leq
J({\bf u})+\lambda_{\mu}\big(J({\bf v})-J({\bf u)}\big)\nonumber \\
= J(\boldsymbol{\mu})+\lambda_{\mu}\big(J(\boldsymbol{\nu})-J(\boldsymbol{\mu})\big).
\end{eqnarray}
In this setting, Algorithm 3.3 uses the convexified cost in Step 3 (Eq. (18)), but by (35), we would get a lower value by
convexifying the control. Consequently, the inequality  in (18) would imply  a similar  an inequality with the convexified control,
as  in Eq. (36), below. Modifying Algorithm 3.3 accordingly, the following algorithm
results.

Given constants $\alpha\in(0,1)$, $\beta\in(0,1)$, and $\eta\in(0,1)$.
\begin{algo}
Given $\boldsymbol{\mu}\thicksim{\bf u}$  such that ${\bf u}$ is an admissible control on $conv(U)$,
compute $\boldsymbol{\mu}_{{\rm next}}\in{\bf M}$ by the following steps.

{\it Step 0:} If $\theta(\boldsymbol{\mu})=0$, set $\boldsymbol{\mu}_{{\rm next}}=\boldsymbol{\mu}$, then exit.\\
{\it Step 1:} Compute the state and costate trajectories, ${\bf x}$ and ${\bf p}$,
associated with $\boldsymbol{\mu}$,
by using Equations (3) and (5), respectively.\\
{\it Step 2:} Compute an $\eta$-minimizer of the Hamiltonian,   $\boldsymbol{\nu}\thicksim{\bf v}$, such that
 ${\bf v}$ is an ordinary control on $conv(U)$. \\
 {\it Step 3:} Compute the integer $\ell_{\mu}$ defined as follows,
 \begin{align}
 & \ell_{\mu}=\min\big\{\ell=0,1,\ldots,: \nonumber \\
 & J({\bf u}+\beta^{\ell}({\bf v}-\bf{u}))-J({\bf u})\leq\alpha
 \beta^{\ell}\eta\theta(\boldsymbol{\mu})\big\}.
 \end{align}
 Define  $\lambda_{\mu}:=\beta^{\ell_{\mu}}$.\\
 {\it Step 4:} Set
 \begin{equation}
  \boldsymbol{\mu}_{{\rm next}}\thicksim{\bf u}+\lambda_{\mu}({\bf v}-{\bf u}).
 \end{equation}
 \end{algo}
 Note that, by (34) and (37),   an iterative application of  Algorithm 3.6
 would compute only ordinary controls that are admissible on  $conv(U)$.
 Furthermore, as discussed earlier, by
 Equation (35), if the Armijo test in Equation (18) were to be satisfied for a given $\ell$ then it would be satisfied in (36) as well
and result in a lower descent in $J(\boldsymbol{\mu}_{{\rm next}})-J(\boldsymbol{\mu})$.
Since this descent in (18) yields the key condition of uniform  descent for Algorithm 3.3, it also
holds for Algorithm 3.6 thereby guaranteeing its convergence via a verbatim application of Proposition 3.4.

\section{Simulation Results}
This section reports on  applications of  Algorithm 3.3 and Algorithm 3.6  to three problems: an autonomous switched-mode
problem, a controlled switched-mode problem, and a problem of balancing motion energy with transmission energy
in a mobile network.\footnote{The
code was written in MATLAB and executed  on
 a laptop computer with  an Intel i7 quad-core processor,   clock frequency of 2.1 GHz, and 8GB of RAM.} In the first problem
the state equation and the cost function are affine in $u$ and hence Algorithm 3.3 is identical to Algorithm 3.6;
in the second problem the cost function is not affine in $u$ and hence we use Algorithm 3.6; and in the third problem the
state equation is affine in $u$ and the cost function in convex in $u$ and hence we can use Algorithm 3.6.
The  first
two problems were considered in \cite{Vasudevan12}, and we use its reported results as benchmarks for our algorithms.
The third problem was  addressed in  \cite{Jaleel13}, but we choose here an initial guess that is farther from the optimum.
 As stated earlier the efficiency of the algorithm depends on the ease with which the pointwise minimizer of the
Hamiltonian can be computed, and for all three problems it will be shown to be  computable via a simple formula. Since the resulting function
${\bf u}^{\star}$ is an ordinary  control, we use  $\boldsymbol{\nu}\thicksim{\bf u}$ in
Step 2.

\subsection{Double Tank System}
Consider a fluid-storage tank with a constant horizontal cross section, where fluid
enters
from the top and discharged through a hole at the bottom. Let $v(t)$ denote the fluid inflow rate from the top,
and let $x(t)$ be the fluid level in the tank. According to Toricelli's law, the state equation of the system is
$\dot{x}(t)=v(t)-\sqrt{x(t)}$.  The  system considered in this subsection is comprised of
two such tanks, one on top of the other, where the fluid input to the upper tank is from a valve-controlled  hose at the top,
and the input to the lower tank consists of the outflow process from the upper tank. Denoting by
$u(t)$ the inflow rate to the upper tank,
and by $x(t):=(x_{1}(t),x_{2}(t))^{\top}$ the fluid levels in the upper tank and lower tank, respectively,
the state equation of the system is
\begin{equation}
\dot{x} = \left(\begin{array}{c}
                                     u - \sqrt{x_1} \\
                                     \sqrt{x_1} - \sqrt{x_2} \end{array}\right);
\end{equation}
we assume the initial condition $x(0)=(2.0,2.0)^{\top}$. The control input $u(t)$ is assumed to be constrained
to the two-point set $U:=\{1.0,2.0\}$, and hence the system can be viewed as an autonomous switched-mode system
whose modes correspond to the two possible values of $u$.   The considered problem is to have the
fluid
level at the lower tank track the value of $3.0$, and accordingly we choose  the cost functional $J$ to be
\begin{equation}
J = 2\int_{0}^{t_{f}} (x_2 - 3)^2\,dt.
\end{equation}
As in \cite{Vasudevan12}, the final time is $t_{f}=10.0$.

By (38), $f(x,u)$ is affine in $u$, and by (39), $L(x,u)$ does not  depend on $u$ and hence can be considered
affine. Therefore Algorithm 3.6  can be run as a special case of
Algorithm 3.3, where $conv(U)=[1,2]$.  Moreover,
with the costate $p:=(p_{1},p_{2})^{\top}\in R^2$, the Hamiltonian has the form
\[
H(x,u,p)=p_{1}u-p_{1}\sqrt{x_{1}}+p_{2}(\sqrt{x_{1}}-\sqrt{x_{2}})+2(x_{2}-3)^2,
\]
whose pointwise minimizer is
\[
u^{\star}=\left\{
\begin{array}{ll}
1, & {\rm if}\ p_{1}\geq 0,\\
2, & {\rm if}\ p_{1}<0;
\end{array}
\right.
\]
if $p_{1}=0$ then $u^{\star}$ can be any point in the interval $[1,2]$.

We ran Algorithm 3.6 starting  from the initial control $u_{1}(t)=1$ $\forall t \in [0, t_f]$, having the cost
$J({\bf u}_{1})=50.5457$.  All numerical integrations were performed by the forward Euler method with $\Delta t=0.01$,
and we approximate ${\bf u}^{\star}$ by its zero-order hold with the sample values $u^{\star}(i\Delta t)$,
$i=0,1,\ldots,$.  We benchmark the results against the reported
run of the algorithm in \cite{Vasudevan12} which, starting from the same initial control, obtained the final cost of
4.829; we reached a similar final cost. In fact,
100 iterations of our algorithm reduced the cost from
 $J(\boldsymbol{\mu}_{1})=50.5457$
to $J(\boldsymbol{\mu}_{100})=4.7440$ in
$2.6700$ seconds of CPU time.

Figure 1 depicts the graph of $J(\boldsymbol{\mu}_{k})$ vs. the iteration count $k=1,\ldots,100$. The graph indicates
a rapid reduction in the cost  from its initial value  until it stabilizes
after about $18$ iterations.
The L-shaped graph is not atypical in applications of descent algorithms with Armijo step sizes, whose strength
lies in its global stability  and large strides towards local solution points at the initial
stages of its runs. As a matter of fact, similar L shaped graphs were obtained from all of the algorithm
runs reported on in this section.

\begin{figure}
\centering
\includegraphics[width=3.2in]{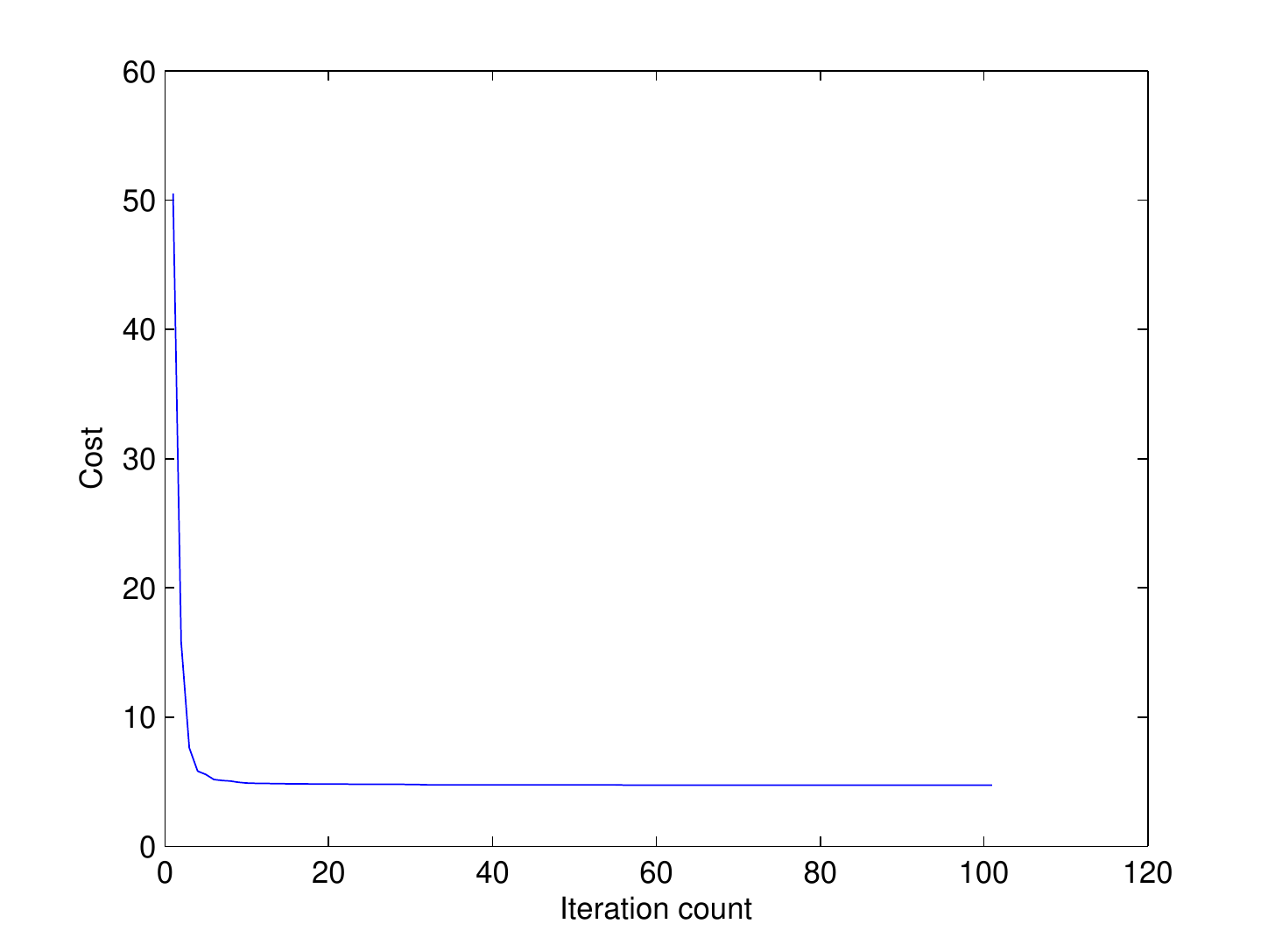}
\caption{ Two-tank system:
$J(\boldsymbol{\mu}_{k})$ vs. $k=1,\ldots,100$; $\Delta t=0.01$}.
\label{fig:two_tank}
\end{figure}

The final relaxed control computed by  Algorithm 3.6,  $\boldsymbol{\mu}_{100}:=\thicksim{\bf u}_{100}$,
was projected onto the space of
admissible (switched-mode)  controls by  Pulse-Width Modulation (PWM) with cycle time of $0.5$ seconds.
% in the following manner.  $u_{100}(t)$
%has the form $u_{100}(t)=\alpha_{1}(t)u_{1}+\alpha_{2}(t)u_{2}$, with $u_{1}:=1$, $u_{2}:=2$; $\alpha_{1}(t)\leq 1$,
%$\alpha_{2}(t)\leq 1$, and $\alpha_{1}(t)+\alpha_{2}(t)=1$.
 %Let us partition the time horizon $[0,t_{f}]$  into
%consecutive $\Delta t$-long cycles. On a typical cycle whose starting time is $\tau$, the projected switched-mode control $u(t)$
%first
%has the value $u(t)=1$ on a subinterval of length $\alpha_{1}(\tau^{+})\Delta t$, and then the value $u(t)=2$ on the remaining,
%$\alpha_{2}(\tau^{+})$-long subinterval.
The  resulting  cost value is  $J({\bf u}_{{\rm fin}})=4.7446$, and the combined runs of the algorithm and the
projection took 2.6825 seconds of CPU time. We
point out that the projection had to be performed only once, after Algorithm 3.6 had completed its
run.
%Alternative, faster and more accurate projection techniques can be found in    \cite{DeCarlo}, but we chose this one due to its simplicity.

Returning to the  run of Algorithm 3.6, the L-shaped graph
in Figure 1 suggests that a reduction in CPU times can be attained, if necessary, by
computing fewer iterations. Moreover, further reduction can be obtained by taking larger integration steps without
significant changes in the final cost. To test this point we ran the algorithm from the
same initial control for 50 iterations with $\Delta t=0.05$,
and it reduced the cost-value from $J(\boldsymbol{\mu}_{1})=50.5282$ to $J(\boldsymbol{\mu}_{50})=4.8078$ in 0.2939 seconds of CPU time;
including the projection onto the space of switched-mode controls it reached $J({\bf u}_{{\rm fin}})=4.8139$
in a total time of 0.3043 seconds.
With a larger integration step,   $\Delta t=0.1$, the algorithm yielded a cost-reduction
from $J(\boldsymbol{\mu}_{1})=50.5069$ to $J(\boldsymbol{\mu}_{50})=4.8816$ in 0.1566 seconds of CPU time, and
$J({\bf u}_{{\rm fin}})=4.8915$ in a total time of 0.1655 seconds. These results are summarized in Table 1.

\begin{table}[h!]
\begin{center}
\begin{tabular}{| c || c | c || c | c |}
\hline
%\multicolumn{3}{|c|}{AAA} & \multicolumn{2}{|c|}{BBB}\\ \hline
$\Delta t;\ k$ & $J(\boldsymbol{\mu}_{k})$ & CPU & $J({\bf u}_{{\rm fin}})$ & CPU\\ \hline \hline
0.01; 100 & 4.7440   & 2.6700   & 4.7446  & 2.6825 \\ \hline
0.05; 50 & 4.8078 & 0.2939 & 4.8139 & 0.3043 \\ \hline
0.1; 50 & 4.8816   & 0.1566   & 4.8915  & 0.1655 \\ \hline
\end{tabular} \par
\bigskip
Table 1: Double-tank problem, $J(\boldsymbol{\mu}_{1})=50.546$
\end{center}
\end{table}

The CPU times indicated in Table 1 are less than the  run-time reported in \cite{Vasudevan12} for  solving
the same problem (32.38 seconds). However, these numbers should not be considered as a sole basis for comparing the
two techniques  since the respective algorithms were implemented on different hardware and software platforms.\footnote{Ref.
\cite{Vasudevan12} refers to a software package for its algorithm, but we were unable to run it because  apparently
it is linked to a proprietary code. Therefore we were unable to conduct a direct comparison between the two algorithms.}
  Furthermore,
the algorithm in  \cite{Vasudevan12} has a broader scope than ours,
while our code is specific for the problem in question. The only conclusion
we draw from Table 1 is that Algorithm 3.6 may have merit and deserves further investigation.
 We believe, however, that  our choice of the
descent direction, namely the pointwise minimizer of the Hamiltonian, plays a role in the fast run times as well as
simplicity of the code as compared with explicit-gradient techniques.

 We close this discussion with a comment on convergence of the algorithm in the control space.
Algorithm 3.6 (as well as Algorithm 3.3) is defined in the space of relaxed controls where its convergence is established
in the weak star topology. Therefore, there is no reason to expect the sequence of computed controls,
$\{{\bf u}_{k}\}$, to converge
 in any (strong) functional norm such as  $L^1$. In fact, the graphs of  $u_{k}(t)$
for $k=1,20,100$, are depicted in Figure 2, where  no strong convergence is discerned. However, the weak
 convergence proved in Section 3 suggests that the cost-sequence $\{J({\bf u}_{k})\}$ would converge
 to the minimal cost, and this indeed is evident from Figure 1. Furthermore, the associated sequence of
 state trajectories are expected to converge in a strong sense ($L^{\infty}$ norm)  to the state trajectory of the optimal
 control,
 and this is indicated by Figure 3 depicting the graphs of $x_{k,2}(t)$ (the fluid levels at the lower tank)
 for $k=1,20,100$.

\begin{figure}
\centering
\includegraphics[width=3.2in]{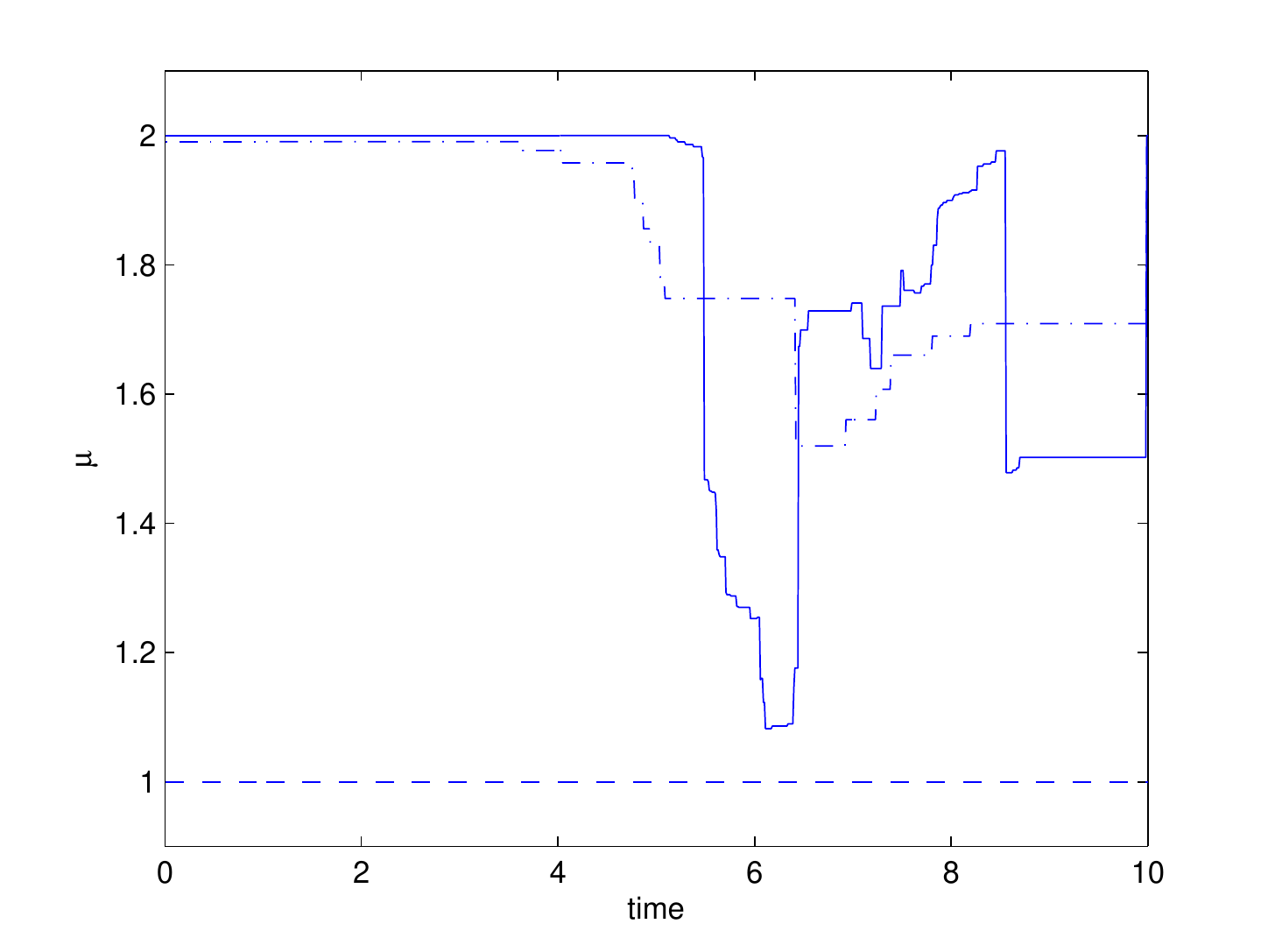}
\caption{Two-tank system:
graph of $u_{1}(t)$ (dashed), $u_{20}(t)$ (dash-dotted), and $u_{100}(t)$ (solid) }
%\label{fig:two_tank}
\end{figure}

\begin{figure}
\centering
\includegraphics[width=3.2in]{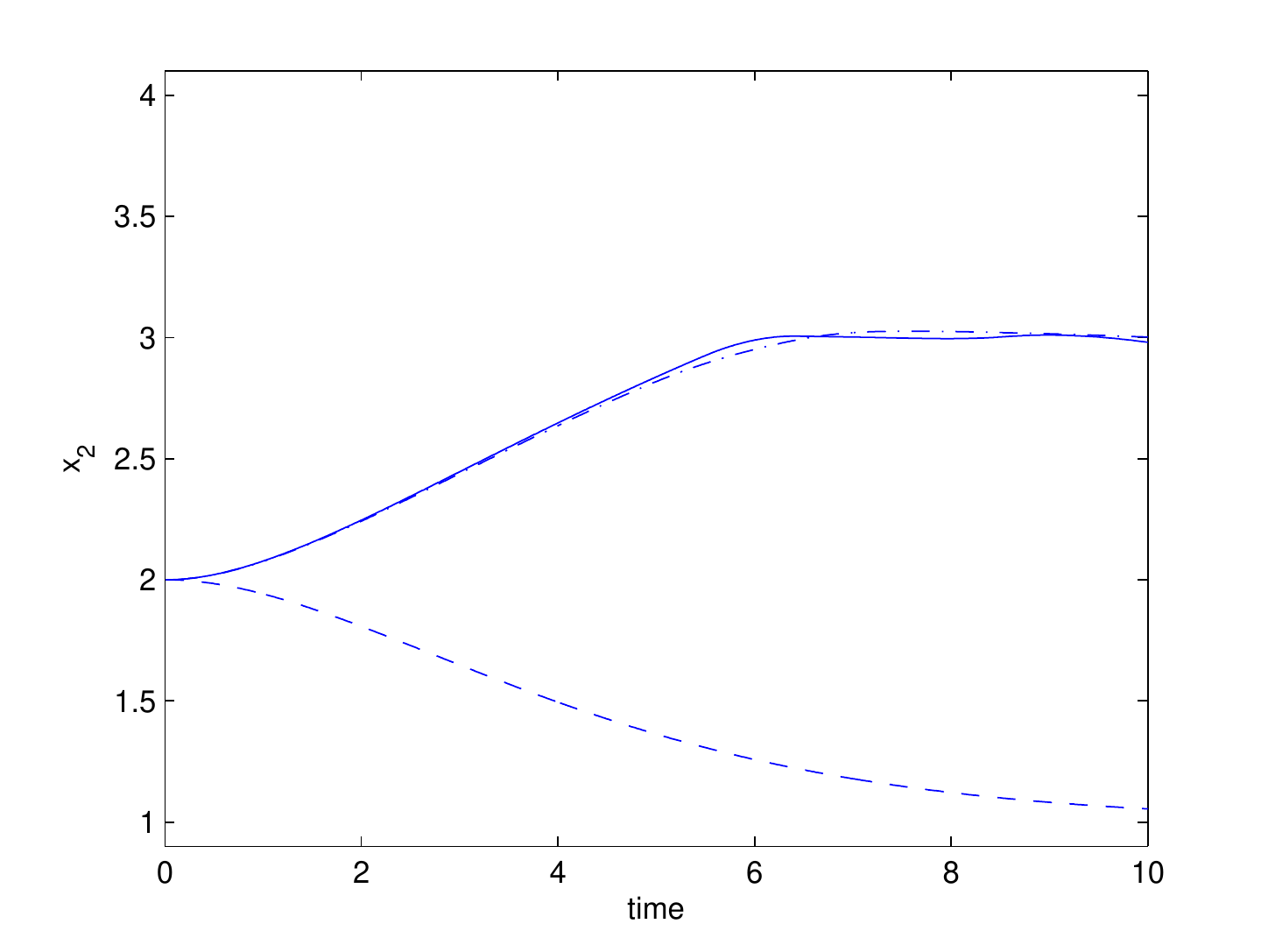}
\caption{Two-tank system:
graph of $x_{2,1}(t)$ (dashed), $x_{2,20}(t)$ (dash-dotted), and $x_{2,100}(t)$ (solid) }
%\label{fig:two_tank}
\end{figure}

\subsection{Hybrid LQR}
This problem was considered in \cite{Vasudevan12} as well. Consider the switched-linear system $\dot{x}=Ax+bv$,
where $x\in R^3$,
\[
A\ =\ \left(\begin{array}{ccc}
                        1.0979 &-0.0105 & 0.0167 \\
                       -0.0105 & 1.0481 & 0.0825 \\
                        0.0167 & 0.0825 & 1.1540 \end{array}\right),
\]
$b\in R^3$ is constrained to a finite  set $B\subset R^3$, and $v\in R$ is a continuum-valued  input. We assume the initial condition to
be
$x_{0}:=x(0)=(0,0,0)^{\top}$.
The set $B$ consists of three points, namely $B=\{b_{1},b_{2},b_{3}\}$, with
$b_1 = (0.9801,-0.1987,0)^{\top}$,  $b_2 = (0.1743,0.8601,-0.4794)^{\top}$, and $b_3 = (0.0952,0.4699,0.8776)^{\top}$. The continuum-valued  control $v$ is constrained to $|v|\leq 20.0$.

All three modes of the system are unstable since the eigenvalues of $A$ are in the right-half plane,  and the problem is to switch among the vectors $b\in B$ and choose $\{v(t)\}$  in a way that brings the state
close to the target state $x_{f}:=(1,1,1)^{\top}$ in a given amount of time and in a way that minimizes the input energy.
The corresponding cost functional is
\begin{equation}
J = 0.01\int_{0}^{t_{f}} v^2\,dt + \|x(t_{f}) - x_{f}\|^2
\end{equation}
with $t_{f}=2.0$, and the control variable is $u=(b,v)\in B\times[-20,20]$.
Reference \cite{Vasudevan12} solved this problem from the initial control of $v_{1}(t)=0$
and $b = (0.9801, -0.1987, 0)^{\top}$ for all $t \in [t_0, t_f]$, and attained
 the final cost of $1.23\cdot 10^{-3}$ in (reported) $9.827$ seconds of CPU time.

Since $f(x,u)$ is not affine in $u$ (due to the term $bv$) Algorithm 3.6 may not be applicable and hence we
used  Algorithm 3.3. We chose the same initial guess  ${\bf u}_{1}$ as in \cite{Vasudevan12} and ran the algorithm for 20 iterations.
The results indicate a similar L-shaped graph of
$J(\boldsymbol{\mu}_{k})$ to the one shown in Figure 1, and attained a cost-value reduction from
$J({\bf u}_{1})=3.000$ to $J(\boldsymbol{\mu}_{20})=
2.768\cdot 10^{-3}$  in $0.761$ seconds of CPU time. All integrations were performed by the forward Euler method
with $\Delta t=0.01$, and the boundary condition of Equation (5) was $p(t_{f})=2\big(x(t_{f})-(1,1,1)^{\top}\big)$ due to the
cost on the final state.\footnote{The cost functional in (40) is not quite in the form of (2) due to the
 addition of the final-state cost term $\phi(x(t_{f})):=||x(t_{f})-x_{f}||^2$.
 However, using the standard transformation  of a Bolza optimal control problem (as in (40)) to a Lagrange problem
 (as in  (2))  and applying the algorithm to the latter, requires only one change,
 namely setting the boundary condition of the costate to $p(t_{f})=\nabla\phi(x_{f})$; all else remains the same.} The Hamiltonian at a control $u=(b,v)$ has the form
$H(x,u,p)=p^{\top}(Ax+bv)+0.01v^2$, and hence its minimizer over $U$, $u^{\star}=(b^{\star},v^{\star})$, is computable as follows:
For every $b_{i}\in B$, $i=1,2,3$, define $v_{i}$ according to the following three contingencies:
(i) if $|p^{\top}b_{i}/0.02|\leq 20$, then $v_{i}=-p^{\top}b_{i}/0.02$;
(ii) if $p^{\top}b_{i}/0.02> 20$, then $v_{i}=-20;$ and (iii) if $p^{\top}b_{i}/0.02< -20$, then $v_{i}=20$.
It is readily seen that $u^{\star}={\rm argmin}\big(H(x,(b_{i},v_{i}),p):i=1,2,3\big)$ is a minimizer of the
Hamiltonian.\footnote{By Proposition 3.2, a pointwise minimizer of the Hamiltonian always can be found in $U$ and does not have to be in $M$.}

 A typical relaxed control can be represented as
$\mu(t)\thicksim\sum_{i=1}^{3}\alpha_{i}(t)b_{i}v_{i}(t)$,
with $\alpha_{i}(t)\in[0,1]$, $i=1,2,3$; $\sum_{i=1}^3\alpha_{i}(t)=1$; and $|v_{i}(t)|\leq 20.0$, $i=1,2,3$
(this is an embedded control as defined in   \cite{Bengea05}).
It can be seen, after some algebra, that this can be represented as
$\mu(t)\thicksim\big(\sum_{i=1}^{3}\gamma_{i}(t)\bar{b}_{i}(t)\big)w(t)$,
with $\gamma_{i}(t)\in[0,1]$, $i=1,2,3$; $\sum_{i=1}^3\gamma_{i}(t)=1$; $\bar{b}_{i}(t)\in\{b_{i},-b_{i}\}$,
 $i=1,2,3$; and $|w(t)|\leq 20.0$. Define $\bar{w}_{i}(t)$ as $\bar{w}_{i}(t)=w$ if $b_{i}(t)=b_{i}$, and
 $\bar{w}_{i}(t)=-w$ if $b_{i}(t)=-b_{i}$.
The projection $\mu(t)$
onto the space of ordinary controls was done by PWM as described in the previous subsection,
with $u$ having the successive values $(b_{1},\bar{w}_{1}(t)),\ (b_{2},\bar{w}_{2}(t)),\ (b_{3},\bar{w}_{3}(t))$ in each cycle
according to the coefficients $\gamma_{i}(t)$, $i=1,2,3$. The cycle time was 12$\Delta t$. Twenty iterations of Algorithm 3.3
followed by the projection of $\boldsymbol{\mu}_{20}$ onto the space of switched-mode controls
required a total CPU time of $0.803$ seconds and yielded a final cost of
$J({\bf u}_{\rm fin})=2.956\cdot 10^{-3}$; the results are summarized in Table 2.

\begin{table}[h!]
\begin{center}
\begin{tabular}{| c || c | c || c | c |}
\hline
%\multicolumn{3}{|c|}{AAA} & \multicolumn{2}{|c|}{BBB}\\ \hline
$\Delta t;\ k$ & $J(\boldsymbol{\mu}_{k})$ & CPU & $J({\bf u}_{{\rm fin}})$ & CPU\\ \hline \hline
0.01; 20 & $2.768\cdot 10^{-3}$   & 0.761   & $2.956\cdot 10^{-3}$  & 0.803  \\ \hline
\end{tabular} \par
\bigskip
Table 2: Hybrid-LQR problem, $J(\boldsymbol{\mu}_{1})=3.00$
\end{center}
\end{table}

Now consider the problem of minimizing the cost functional $0.01\int_{0}^{t_f}v^2dt$ subject to the constraint
$x(t_f)=x_f:=(1.0,1.0,1.0)^{\top}$. The definition of $J$ in Eq. (40) appears to address this problem with the penalty function $||x(t_f)-x_f||^2$. As a matter of fact,
the final state $x(t_f)$ obtained form the run of the algorithm is $x(t_f)=(0.9994,0.9991,0.9998)^{\top}$,
and after projecting the  relaxed control ${\bf u}_{20}$ onto the space of ordinary controls, the corresponding final state is
$(0.9965,1.0020,0.9875)^{\top}$.

While Algorithm 3.3 could be applied to the current problem, its scope does not include some embedded optimal control problems, comprising a class
or relaxed-control problems defined in \cite{Bengea05}. However, we believe that a possible  extension of Algorithm 3.6 can close this gap, and is
currently under investigation.

\subsection{Balancing Mobility with Transmission Energy  in Mobile Sensor Networks}
In Reference \cite{Jaleel13} we considered a path-planning problem for mobile communication-relay networks, whose objective is
to optimize a weighted sum of transmission energy and fuel consumption.
We used there Algorithm 3.6 to solve it. Next we present simulation results for the same
problem, but with a different initial control, ${\bf u}_{1}$, chosen farther from the optimum in order to highlight
 the drastic cost-reduction of the algorithm's run  at its initial phases.

Consider a scenario where
a given number ($N$) of mobile sensors (agents) are placed in a terrain, and at time $t=0$ they are tasked with
forming a point-to-point relay network for communications between a stationary object and a stationary controller. Upon issuance of the command
the agents start moving while transmitting. Given the final time $t_{f}$, the  problem is to compute the agents' paths in a way
that minimizes a weighted sum of their fuel consumption  and transmission energy over the time-interval
$t\in[0,t_{f}]$. Of course the optimal paths depend on the positions of the object and controller, as well as on
the initial positions of the agents.

A detailed description of the problem and justification of the assumptions made can be found
in \cite{Jaleel13}. As in \cite{Jaleel13},
we assume  that the agents' positions $x_{i}$, $i=1,\ldots,N$, are confined to a line-segment
$[0,d]$ for a given $d>0$.
 Denote by   $x=(x_{1},\ldots,x_{N})^{\top}\in R^N$ the vector of the agents' positions,
 and by $u:=(u_{1},\ldots,u_{N})^{\top}\in R^N$,  the vector of their corresponding   velocities.
Viewing $x$ as the state of the system and $u$ as its control input, the state equation is
\begin{equation}
\dot{x}=u.
\end{equation}
The power required to transmit a signal over a $z$-long channel can be considered as
proportional to $z^2$ (see \cite{Jaleel13}), and the fuel rate required to move an agent is proportional to the speed of
motion.  Consequently, and defining $x_{0}:=0$ and $x_{N+1}:=d$,    the considered cost-performance functional
is
\begin{equation}
J\ =\ \sum_{i=1}^{N+1}\int_{0}^{t_{f}}(x_{i}-x_{i-1})^2dt+C\sum_{i=1}^{N}\int_{0}^{t_{f}}|u_{i}|dt
\end{equation}
for a given $C>0$. The optimal control problem is to minimize $J$ for a given initial condition $x(0)$, subject to the
pointwise input constraints $|u_{i}|\leq\bar{u}$ for a given $\bar{u}>0$.

By Equations (41)-(42) it can be seen that the costate $p:=(p_{1},\ldots,p_{N})^{\top}$ is defined by the  equation
\begin{equation}
\dot{p}_{i}=2(x_{i-1}+x_{i+1}-2x_{i}),
\end{equation}
$i=1,\ldots,N$, with the boundary condition
$p_{i}(t_{f})=0$.
Therefore,
given a control $u\in R^N$ and its associated state $x\in R^N$ and costate
$p=(p_{1},\ldots,p_{N})^{\top}\in R^N$, the Hamiltonian has the
form
$H(x,u,p)=\sum_{i=1}^Np_{i}u_{i}+J$ with $J$ defined in (42), and its minimizer,  $u^{\star}=(u_{1}^{\star},\ldots,u_{N}^{\star})^{\top}$, is computable as follows,
\begin{equation}
u_{i}^{\star}=\left\{
\begin{array}{ll}
-{\rm sgn}(p_{i})\bar{u}, & {\rm if}\ |p_{i}|>C\\
0, & {\rm if}\ |p_{i}|\leq C
\end{array}
\right.
\end{equation}
(e.g., \cite{Jaleel13}).

 In our simulation experiments  we considered an example with $N=6$ (six agents), $d=20$ (hence the agents move in the interval $[0,20]$),
 $t_{f}=20$, $C=7$, $\bar{u}=1$, and the initial
 state is $x(0)=(1,2,7,9,12,19)^{\top}$. The numerical integrations were performed by the forward Euler method with
 the step size $\Delta t=0.01$. The algorithm started with the following initial control,
 ${\bf u}_{1}:=({\bf u}_{1,1},\ldots,{\bf u}_{1,6})^{\top}$:
$u_{1}(t)=1.0$;
$u_{2}(t)=sin(\pi t/4)$;
$u_{3}(t)=3u_{2}(t)$;
$u_{4}(t)=2u_{3}(t)$;
$u_{5}(t)=2u_{4}(t)$; and
$u_{6}(t)=u_{5}(t)-4.3$, with the initial cost of
$J(u_{1})=81,883.4$.

A 200-iteration run took 11.3647 seconds of CPU time and yielded
  the final cost of  $J({\bf u}_{200})=1,253.4$
  (extensive simulations in \cite{Jaleel13} suggested that this is about the global minimum).  The graph of $J({\bf u}_{k})$ vs. $k$ has a similar L shape to Figure 1, and it took only 4 iterations (9 iterations, resp.) to achieve 98\% (99\%, resp.) of the total cost reduction to $J({\bf u}_{5})=2,701.6$ ($J({\bf u}_{10})=2,037.6$, resp).
  To reduce the CPU times  we can take fewer iterations.  For example, 20 iterations take 0.8653 seconds of CPU times to obtain $J({\bf u}_{20})=1,455.5$, and 100 iterations took 5.2379 seconds to yield $J({\bf u}_{100})=1,256.7$.
  Further speedup can be achieved by increasing the integration step size: with $\Delta t=0.1$, 100
  iterations took 0.5645 seconds of CPU time to obtain $J({\bf u}_{100})=1,260.4$. Note that this is quite close to the
  aforementioned apparent minimum of 1,253.4. These results are summarized in Table 3.\\

\begin{table}[h!]
\begin{center}
\begin{tabular}{| c || c | c  | c |}
\hline
$\Delta t$ & $\ k$ & $J({\bf\mu}_{k})$ & CPU \\ \hline \hline
0.01 & 200 & 1,253.4 & 11.3647 \\ \hline
0.01 & 100 & 1,256.7 & 5.2379 \\ \hline
0.01 & 20 & 1,455.5 & 0.8653 \\ \hline
0.1 & 100 & 1,260.4 & 0.5645 \\ \hline
\end{tabular} \par
\bigskip
Table 3: Path planning for power-aware mobile networks, $J({\bf u}_{1})=81,883.4$
\end{center}
\end{table}

  \section{Conclusions}
This paper presents an iterative algorithm for solving a class of optimal control problems. The algorithm operates in the
space of relaxed controls and the obtained  result is projected onto the space of ordinary
controls. The computation of the descent direction is based on pointwise minimization of the Hamiltonian at each iteration instead of explicit gradient calculations. Simulation examples indicate fast convergence for a number of test problems.

\bibliography{biblio}

\begin{thebibliography}{14}
\providecommand{\natexlab}[1]{#1}
\providecommand{\url}[1]{\texttt{#1}}
\providecommand{\urlprefix}{URL }
\expandafter\ifx\csname urlstyle\endcsname\relax
  \providecommand{\doi}[1]{doi:\discretionary{}{}{}#1}\else
  \providecommand{\doi}{doi:\discretionary{}{}{}\begingroup
  \urlstyle{rm}\Url}\fi



\bibitem{Polak97}
E. Polak. {\it Optimization Algorithms and Consistent
Approximations}. Springer-Verlag, New York, New York, 1997.


\bibitem{Caldwell10}
T. Caldwell and T. Murphy.
An Adjoint Method for Second-Order Switching Time Optimization.
{\it Proc. 49th CDC}, Atlanta, Georgia, December 15-17, 2010.


\bibitem{Caldwell11}
T. Caldwell and T. D. Murphey. Switching mode generation and optimal
estimation with application to skid-steering, {\it Automatica}, vol. 47,
no. 1, pp. 50–64, 2011.

\bibitem{Caldwell12}
T. Caldwell and T. D. Murphey. Single integration optimization of
linear time-varying switched systems. {\it IEEE Transactions on Automatic
Control}, vol. 57, no. 6, pp. 1592–1597, 2012.

\bibitem{Caldwell16}
T. Caldwell and T. Murphey.  Projection-Based Iterative Mode
Scheduling for Switched Systems, {\it Nonlinear Analysis: Hybrid Systems},
to appear, 2016.


\bibitem{Gonzalez10}
H. Gonzalez, R. Vasudevan, M. Kamgarpour, S.S. Sastry, R. Bajcsy, and C. Tomlin.
A Numerical Method for the Optimal Control of Switched Systems.
{\it Proc. 49th CDC}, Atlanta, Georgia, pp. 7519-7526, December 15-17, 2010.

\bibitem{Vasudevan12}
R. Vasudevan, H. Gonzalez, R. Bajcsy, and S.S. Sastry. Consistent
Approximations for the Optimal Control of Constrained Switched
Systems - Part 1: A Conceptual Algorithm, and Part 2: An Implementable Algorithm.  {\it SIAM Journal on Control and Optimization}, Vol. 51, pp.
4663-4483 (Part 1) and pp. 4484-4503 (Part 2), 2013.


\bibitem{Egerstedt06}
M. Egerstedt, Y. Wardi, and H. Axelsson. Transition-Time
Optimization for Switched Systems.  {\it IEEE Transactions on
Automatic Control}, Vol. AC-51, No. 1, pp. 110-115, 2006.

\bibitem{Axelsson08}
H. Axelsson, Y. Wardi, M. Egerstedt, and E. Verriest.
 A Gradient Descent Approach to Optimal Mode Scheduling in Hybrid Dynamical Systems.
 {\it  Journal of Optimization Theory and Applications}, Vol. 136,  pp. 167-186,  2008.

\bibitem{McShane67}
E.J. McShane.
Ralexed Controls and Variational Problems.
{\it SIAM Journal on Control},
Vol. 5, pp. 438-485, 1967.

\bibitem{Warga72}
J. Warga,
{\it Optimal Control of Differential and Functional Equations},
Academic Press, 1972.

\bibitem{Armijo66}
L. Armijo. Minimization of Functions Having Lipschitz Continuous
First-Partial Derivatives.
{\it Pacific Journal of Mathematics}, Vol. 16, pp. 1-3, 1966.

\bibitem{Branicky98}
M.S. Branicky, V.S. Borkar, and S.K. Mitter. A Unified Framework
for Hybrid Control: Model and Optimal Control Theory. {\it IEEE
Transactions on Automatic Control}, Vol. 43, pp. 31-45, 1998.

\bibitem{Piccoli98}
B. Piccoli. Hybrid Systems and Optimal Control. {\it Proc. IEEE
Conference on Decision and Control}, Tampa, Florida, pp. 13-18,
1998.

\bibitem{Sussmann99}
H.J. Sussmann. A Maximum Principle for Hybrid Optimal Control
Problems. {\it Proceedings of the 38th IEEE Conference on Decision
and Control}, pp. 425-430, Phoenix, AZ, Dec. 1999.


\bibitem{Shaikh02}
M.S. Shaikh and P. Caines. On Trajectory Optimization for Hybrid
Systems: Theory and Algorithms for Fixed Schedules. {\it IEEE
Conference on Decision and Control}, Las Vegas, NV, Dec. 2002.

\bibitem{Garavello05}
M. Garavello and  B. Piccoli.
 Hybrid Necessary Principle. {\it  SIAM J. Control and Optimization}, Vol.  43, pp. 1867-1887,
2005.



 \bibitem{Bengea05}
S.C. Bengea and R. A. DeCarlo.
Optimal control of switching systems.
{\it Automatica}, Vol. 41, pp. 11-27, 2005.

\bibitem{Shaikh07}
M.S. Shaikh and P.E. Caines.
On the Hybrid Optimal Control Problem: Theory and Algorithms.
{\it IEEE Trans. Automatic Control}, Vol. 52, pp. 1587-1603, 2007.


\bibitem{Taringoo13}
F. Taringoo and P.E.  Caines.
On the optimal control of hybrid systems on Lie groups and the exponential gradient HMP algorithm.
In {\it Proc. 52nd IEEE Conf. on Decision and Control}, Florence, Italy, December 10-13, 2013.

\bibitem{Taringoo13a}
F. Taringoo and P.E.  Caines.
On the Optimal Control of Impulsive Hybrid Systems on Riemannian Manifolds.
{\it SIAM Journal on Control and Optimization},  Vol. 51, Issue 4, pp. 3127 - 3153, 2013.



\bibitem{Xu02}
X. Xu and P. Antsaklis. Optimal Control of Switched Autonomous
Systems. {\it IEEE Conference on Decision and Control}, Las Vegas,
NV, Dec. 2002.

\bibitem{Xu02a}
X. Xu and P.J. Antsaklis. Optimal Control of Switched Systems via
Nonlinear Optimization Based on Direct Differentiations of Value
Functions. {\it International Journal of Control}, Vol. 75, pp.
1406-1426, 2002.


\bibitem{Caldwell12a}
T. M. Caldwell and T. D. Murphey.
Projection-based Switched System Optimization.
{\it Proc. American Control Conference}, Montreal, Canada, June 2012.

\bibitem{Caldwell12b}
T. M. Caldwell and T. D. Murphey. Projection-Based Switched
System Optimization: Absolute Continuity of the Line Search.  {\it Proc.
51st CDC}, Maui, Hawaii, December 10-13, 3012.









\bibitem{Shaikh03}
M.S. Shaikh and P.E. Caines. On the Optimal Control of Hybrid
Systems: Optimization of Trajectories, Switching Times and
Location Schedules. In {\it Proceedings of the 6th International
Workshop on Hybrid Systems: Computation and Control}, pp. 466-481,
Prague, The Czech Republic, 2003.

\bibitem{Caines05}
P. Caines and M.S. Shaikh. Optimality Zone Algorithms for Hybrid
Systems Computation and Control: Exponential to Linear Complexity.
{\it Proc. 13th Mediterranean Conference on Control and
Automation}, Limassol, Cyprus, pp. 1292-1297, June 27-29, 2005.

\bibitem{Shaikh05}
M.S. Shaikh and P.E. Caines.
Optimality Zone Algorithms for Hybrid Systems Computation and Control:
From Exponential to Linear Complexity.
{\it Proc. IEEE Conference on Decision and Control/European Control
Conference}, pp. 1403-1408, Seville, Spain, December 2005.

\bibitem{Hellund99}
S. Hedlund and A. Rantzer.
Optimal control for hybrid systems.
{\it Proc. 38th CDC},
Phoenix, Arizona, December 7-10, 1999.

\bibitem{Attia05}
S.A. Attia, M. Alamir, and C. Canudas de Wit. Sub Optimal Control of
Switched Nonlinear Systems Under Location and Switching Constraints.
{\it Proc. 16th IFAC World Congress}, Prague, the Czech Republic,
July 3-8, 2005.


\bibitem{Ge75}
X. Ge, W. Kohn, A. Nerode, and J.B. Remmel.
Hybrid systems: Chattering approximation to
relaxed controls.
{\it Hybrid Systems III: Lecture Notes in Computer Science},
R. Alur, T. Henzinger, E. Sontag, eds., Springer Verlag,
Vol. 1066, pp. 76-100, 1996.

%\bibitem{DeCarlo}
%R.T. Meyer, R.A.  DeCarlo, P.H.  Meckl,  and S. Pekarek.
%Hybrid Model
%Predictive Power Management of a Fuel-cell Battery Vehicle.
%{\it  Asian Journal of
%Control}, Vol.  15, no. 3, pp. 1-17,  2013.


\bibitem{Meyer12}
R.T. Meyer, M. Zefran, and R.A. Decarlo. Comparison of the Embedding Method to
Multi-Parametric Programming, Mixed-Integer Programming, Gradient Descent, and
Hybrid Minimum Principle Based Methods.
{\it IEEE Transactions on Control Systems
Technology},
Vol. 22, no. 5,  pp. 1784-1800, 2014.


\bibitem{Zhu11}
F. Zhu and P.J. Antsaklis. Optimal Control of Switched Hybrid Systems:
A Brief Survey. Technical Report of the ISIS Group at the University
of Notre Dame. ISIS-2011-003, July 2011.


\bibitem{Lin14}
H. Lin and P. J. Antsaklis.
Hybrid Dynamical Systems:
An Introduction to Control and Verification.
{\it Foundations and Trends  in Systems and Control},
Vol. 1, no. 1, pp. 1-172, March 2014.



%\bibitem{Wardi12}
%Y. Wardi and M. Egerstedt.
%Algorithm for Optimal Mode Scheduling in  Switched Systems.
%{\it Proceedings of the American Control Conference}, Montreal, Canada, June 2012.


\bibitem{Young69}
 L. C. Young.
{\it Lectures on the calculus of variations and optimal control theory}.
Foreword by Wendell H. Fleming. W. B. Saunders Co., Philadelphia, 1969.


\bibitem{Gamkrelidze78}
R. Gamkrelidze. {\it Principle of Optimal Control Theory.} Plenum, New York, 1978.



\bibitem{Berkovitz13}
L.D. Berkovitz and N.G. Medhin.
{\it Nonlinear Optimal Control Theory},
Chapman \& Hall, CRC Press,
Boca Raton, Florida, 2013.

\bibitem{Lou09}
H. Lou.
Analysis of the Optimal Relaxed Control to an Optimal Control Problem.
{\it Applied Mathematics and Optimization}, Vol. 59, pp. 75-97, 2009.

\bibitem{Vinter00}
R. Vinter.
{\it Optimal Control}, Birkhauser, Boston, Massachusetts, 2000.

\bibitem{Mayne75}
D.Q. Mayne and E. Polak.
First-order Strong Variation Algorithms for Optimal Control.
{\it J. Optimization Theory and Applications},
Vol. 16,
pp. 277-301, 1975.

\bibitem{Jaleel13}
H. Jaleel, Y. Wardi, and M. Egerstedt.
Minimizing Mobility and Communication Energy in Robotic Networks: an Optimal Control Approach.
{\it Proc. 2014 American Control Conference}, Portland, Oregon, June 4-6, 2014.


















%\bibitem{Riedinger99}
%P. Riedinger, F. Kratz, C. Iung, and C. Zanne. Linear Quadratic
%Optimization for Hybrid Systems. {\it Proceedings of the 38th IEEE
%Conference on Decision and Control}, pp. 3059-3064, Phoenix,
%Arizona, pp. 3059-3064, December 1999.






















\end{thebibliography}

\bibliographystyle{plain}

\end{document}